\documentclass{article}
\pdfoutput=1

\usepackage{arxiv}

\usepackage[utf8]{inputenc} 
\usepackage[T1]{fontenc}    
\usepackage{hyperref}       
\usepackage{url}            
\usepackage{booktabs}       
\usepackage{amsfonts}       
\usepackage{nicefrac}       
\usepackage{microtype}      
\usepackage{lipsum}		    
\usepackage{amsmath}
\usepackage{times}
\usepackage{hyperref}
\usepackage{float}
\usepackage{enumitem}
\usepackage{comment}

\RequirePackage[numbers]{natbib}
\usepackage{graphicx}

\usepackage{todonotes}

\title{Fair Coins Tend to Land on the Same Side They Started: Evidence from $350{,}757$ Flips} 

\author{
\textbf{František Bartoš$^{1*}$, Alexandra Sarafoglou$^{1}$, Henrik R. Godmann$^{1}$, Amir Sahrani$^{1}$, David Klein Leunk$^{2}$,}\\
\textbf{Pierre Y. Gui$^{2}$, David Voss$^{2}$, Kaleem Ullah$^{2}$, Malte J. Zoubek$^{3}$, Franziska Nippold,}\\
\textbf{Frederik Aust$^{1}$, Felipe F. Vieira$^{4}$, Chris-Gabriel Islam$^{5,6}$, Anton J. Zoubek$^{7}$, Sara Shabani$^8$,}\\
\textbf{Jonas Petter$^{1}$, Ingeborg B. Roos$^9$, Adam Finnemann$^{1,10}$, Aaron B. Lob$^{2,11}$, Madlen F. Hoffstadt$^{1}$,}\\
\textbf{Jason Nak$^{2}$, Jill de Ron$^{2}$, Koen Derks$^{12}$, Karoline Huth$^{2,13}$, Sjoerd Terpstra$^{14}$,}\\
\textbf{Thomas Bastelica$^{15,16}$, Magda Matetovici$^{2,17}$, Vincent L. Ott$^{2}$, Andreea S. Zetea$^{2}$, Katharina Karnbach$^{2}$,}\\
\textbf{Michelle C. Donzallaz$^{1}$, Arne John$^{2}$, Roy M. Moore$^{2}$, Franziska Assion$^{1,8}$, Riet van Bork$^{2}$,}\\
\textbf{Theresa E. Leidinger$^{2}$, Xiaochang Zhao$^{2}$, Adrian Karami Motaghi$^{2}$, Ting Pan$^{1,9}$, Hannah Armstrong$^{2}$,}\\
\textbf{Tianqi Peng$^{2}$, Mara Bialas$^{20}$, Joyce Y.-C. Pang$^{2}$, Bohan Fu$^{2}$, Shujun Yang$^{2}$,}\\
\textbf{Xiaoyi Lin$^{2}$, Dana Sleiffer$^{2}$, Miklos Bognar$^{21}$, Balazs Aczel$^{22}$, and Eric-Jan Wagenmakers$^{1}$}\\
\\
\footnotesize{$^{1}$Department of Psychological Methods, University of Amsterdam,}\\ 
\footnotesize{$^{2}$Department of Psychology, University of Amsterdam,}\\
\footnotesize{$^{3}$Institute of Economic Law, University of Kassel,}\\
\footnotesize{$^{4}$Faculty of Psychology and Educational Sciences, KU Leuven,}\\
\footnotesize{$^{5}$Chairs of Statistics and Econometrics, Georg August University of Göttingen,}\\
\footnotesize{$^{6}$Centre for Environmental Sciences, Hasselt University,}\\
\footnotesize{$^{7}$Gutenberg School Wiesbaden,}\\
\footnotesize{$^{8}$Department of Psychology, Justus Liebig University,}\\
\footnotesize{$^{9}$Department of Adult Psychiatry, Amsterdam University Medical Centre,}\\
\footnotesize{$^{10}$Centre for Urban Mental Health, University of Amsterdam,}\\
\footnotesize{$^{11}$Cognitive and Behavioral Decision Research, University of Zurich,}\\
\footnotesize{$^{12}$Center of Accounting, Auditing \& Control, Nyenrode Business University,}\\
\footnotesize{$^{13}$Department of Psychiatry, Amsterdam UMC, University of Amsterdam,}\\
\footnotesize{$^{14}$Institute for Marine and Atmospheric Research Utrecht, Utrecht University,}\\
\footnotesize{$^{15}$Centre Hospitalier Le Vinatier,}\\
\footnotesize{$^{16}$INSERM, U1028; CNRS, UMR5292, Lyon Neuroscience Research Center,}\\
\footnotesize{Psychiatric Disorders: from Resistance to Response Team,}\\
\footnotesize{$^{17}$Research Institute of Child Development and Education, University of Amsterdam,}\\
\footnotesize{$^{18}$Faculty of Behavioral and Movement Sciences, Vrije Universiteit Amsterdam,}\\
\footnotesize{$^{19}$Department of Develomental Psychology, University of Amsterdam,}\\
\footnotesize{$^{20}$Behavioural Science Institute, Radboud University,}\\
\footnotesize{$^{21}$Doctoral School of Psychology, ELTE Eotvos Lorand University,}\\
\footnotesize{$^{22}$Institute of Psychology, ELTE Eotvos Lorand University}\\
\\
\footnotesize{*Correspondence concerning this article should be addressed to František Bartoš at f.bartos96@gmail.com.}\\
\footnotesize{For all but the last four positions, the authorship order aligns with the number of coin flips contributed.}\\
}

\date{}

\renewcommand{\shorttitle}{$350{,}757$ Coin Flips}

\begin{document} 

\maketitle 
\newpage

\begin{abstract}
Many people have flipped coins but few have stopped to ponder the statistical and physical intricacies of the process. We collected $350{,}757$ coin flips to test the counterintuitive prediction from a physics model of human coin tossing developed by Diaconis, Holmes, and Montgomery (DHM; 2007). The model asserts that when people flip an ordinary coin, it tends to land on the same side it started---DHM estimated the probability of a same-side outcome to be about 51\%. Our data lend strong support to this precise prediction: the coins landed on the same side more often than not, $\text{Pr}(\text{same side}) = 0.508$, 95\% credible interval (CI) [$0.506$, $0.509$], $\text{BF}_{\text{same-side bias}} = 2359$. Furthermore, the data revealed considerable between-people variation in the degree of this same-side bias. Our data also confirmed the generic prediction that when people flip an ordinary coin---with the initial side-up randomly determined---it is equally likely to land heads or tails: $\text{Pr}(\text{heads}) = 0.500$, 95\% CI [$0.498$, $0.502$], $\text{BF}_{\text{heads-tails bias}} = 0.182$. Furthermore, this lack of heads-tails bias does not appear to vary across coins. Additional analyses revealed that the within-people same-side bias decreased as more coins were flipped, an effect that is consistent with the possibility that practice makes people flip coins in a less wobbly fashion. Our data therefore provide strong evidence that when some (but not all) people flip a fair coin, it tends to land on the same side it started.
\end{abstract}

\section*{Introduction}

A coin flip---the act of spinning a coin into the air with your thumb and then catching it in your hand---is often considered the epitome of a chance event. It features as a ubiquitous example in textbooks on probability theory and statistics \citep{kerrich1946experimental, feller1957introduction, jaynes2003probability, kuchenhoff2008coin, gelman2002you} and constituted a game of chance (`capita aut navia'--`heads or ships') already in Roman times (Macrobius \cite{macrobius431saturnalia}, $\sim 431$ AD, 1.7:22).

The simplicity and perceived fairness of a coin flip, coupled with the widespread availability of coins, may explain why it is often used to make even high-stakes decisions. For example, a coin flip was used to determine which of the Wright brothers would attempt the first flight in 1903; who would get the last plane seat for the tour of rock star Buddy Holly (which crashed and left no survivors) in 1959;  who would be the winner of the 1968 European Championship semi-final soccer match between Italy and the Soviet Union (an event which Italy went on to win); which of two companies would be awarded a public project in Toronto in 2003; to break the tie in local political elections in the Philippines in both 2004 and 2013; finally, this year a coin flip was used to determine whether a \$5,000,000 game show prize would be doubled or whether the contestant would be eliminated.

Despite the widespread popularity of coin flipping, few people pause to reflect on the notion that the outcome of a coin flip is anything but random: a coin flip obeys the laws of Newtonian physics in a relatively transparent manner \citep{jaynes2003probability}. According to the standard model of coin flipping \citep{vulovic1986randomness, keller1986probability, engel1992road, strzalko2008dynamics, strzalko2010understanding}, the flip is a deterministic process and the perceived randomness originates from small fluctuations in the initial conditions (regarding starting position, configuration, upward force, and angular momentum) combined with narrow boundaries on the outcome space. Therefore the standard model predicts that when people flip a fair coin, the probability of it landing heads is 50\% (i.e., there is no `heads-tails bias'; conversely, if a coin would land on one side more often than the other, we would say there is a `heads-tails bias').\footnote{Some even assert that a biased coin is a statistical unicorn---everyone talks about it but no one has actually encountered one \citep{gelman2002you}. Physics models support this assertion as long as the coin is not bent \citep{woo2010bent-b} or allowed to spin on the ground \citep{jaynes2003probability, kuchenhoff2008coin}.} 




The standard model of coin flipping was extended by Diaconis, Holmes, and Montgomery (DHM; \cite{diaconis2007dynamical}) who proposed that when people flip a ordinary coin, they introduce a small degree of `precession' or wobble---a change in the direction of the axis of rotation throughout the coin's trajectory. According to the DHM model, precession causes the coin to spend more time in the air with the initial side facing up. Consequently, the coin has a higher chance of landing on the same side as it started (i.e., `same-side bias'). Under the DHM model, this same-side bias is absent only when there is no wobble whatsoever, as any non-zero angle of rotation results in a same-side bias (with a higher degree of wobble resulting in a more pronounced bias). Based on a modest number of empirical observations (featuring coins with ribbons attached and high-frame-rate video recordings) Diaconis et al. \cite{diaconis2007dynamical} measured the off-axis rotations in typical human flips. Based on these observations, the DHM model predicted that a coin flip should land on the same side as it started with a probability of approximately 51\%, just a fraction higher than chance.

Throughout history, several researchers have collected thousands of coin flips. In the 18\textsuperscript{th} century, the famed naturalist Count de Buffon \cite{Buffon1777Essai} collected $2{,}048$ uninterrupted sequences of `heads' in what is possibly the first statistical experiment ever conducted. In the 19\textsuperscript{th} century, the statistician Karl Pearson \cite{pearson1897the} flipped a coin $24{,}000$ times to obtain $12{,}012$ tails. And in the 20\textsuperscript{th} century, the mathematician John Kerrich \cite{kerrich1946experimental} flipped a coin $10{,}000$ times for a total of $5{,}067$ heads. These experiments do not allow a test of the DHM model, however, mostly because it was not recorded whether the coin landed on the same side that it started. A notable exception is a sequence of $40{,}000$ coin flips collected by Janet Larwood and Priscilla Ku \citep{noauthor_40000_nodate}: Larwood always started the flips heads-up, and Ku always tails-up. Whereas Larwood's $10{,}231$ out of $20{,}000$ heads-to-heads flips are indicative of a same-side bias, Ku's $10{,}014$ out of $20{,}000$ tails-to-tails flips do not. Taken together, the $40{,}000$ coin tosses were interpreted to ``yield ambiguous evidence for dynamical bias'' \citep{noauthor_40000_nodate}.


In order to carry out a diagnostic empirical test of the same-side bias hypothesized by DHM, we collected a total of $350{,}757$ coin flips, a number that---to the best of our knowledge---dwarfs all previous efforts. To anticipate our main results, the data reveal overwhelming statistical evidence for the presence of same-side bias. However, this effect needs to be qualified in the sense that it varies across individuals. Moreover, the bias appears to wane with practice. The data also yield moderate evidence for the complete absence of a heads-tails bias. The appendices demonstrate that the same conclusions obtain under a wide range of alternative analysis strategies.

\section*{Methods}
\label{app:methods}


We collected data in three different settings using the same standardized protocol. First, a group of five bachelor students collected at least $15{,}000$ coin flips each as a part of their bachelor thesis project, contributing $75{,}036$ coin flips in total. Second, we organized a series of on-site ``coin flipping marathons'' where 35 people spent up to 12 hours coin-flipping (see e.g., [blinded for review] for a video recording of one of the events), contributing a total of $203{,}440$ coin flips.\footnote{Including $2{,}700$ coin flips collected by the first two authors on a separate occasion.} Third, we issued a call for collaboration via Twitter, which resulted in an additional seven people contributing a total of $72{,}281$ coin flips. 

We encouraged people to flip coins of various currencies and denomination to ascertain the generalizability of the effect. Furthermore, we encouraged coin tossers to exchange coins, as this potentially allows people-specific effects to be disentangled from coin-specific effects. Overall, a group of 48 people\footnote{One of the bachelor students collected data with a family member who is counted as a ``coin-tosser'' but who did not qualify for co-authorship.} (i.e., all but three of the co-authors) tossed coins of 44 different currencies $\times$ denominations and obtained a total number of $350{,}757$ coin flips.

The protocol required that each person collects sequences of 100 consecutive coin flips.\footnote{Some sequences slightly varied in length due to issues with keeping track of the number of flips.} In each sequence, people randomly (or according to an algorithm) selected a starting position (heads-up or tails-up) of the first coin flip, flipped the coin, caught it in their hand, recorded the landing position of the coin (heads-up or tails-up), and proceeded with flipping the coin starting from the same side it landed in the previous trial (we decided for this ``autocorrelated'' procedure as it simplified recording of the outcomes). In case the coin was not caught in hand, the flip was designated as a failure, and repeated from the same starting position. To simplify the recording and minimize coding errors, participants usually marked sides of the coins with permanent marker. To safeguard the integrity of the data collection effort, all participants videotaped and uploaded recordings of their coin flipping sequences.\footnote{There are occasional missing recordings due to failures of recording apparatus/lost files.} See [blinded for review] for the data and video recordings.

\section*{Analysis}

\subsection*{Same-side bias}
An initial analysis confirms the prediction from the DHM model: the coins landed how they started more often than 50\%. Specifically, the data feature $178{,}079$ same-side landings out of $350{,}757$ tosses, $\text{Pr}(\text{same side}) = 0.5077$, 95\% central credible interval (CI) [$0.5060$, $0.5094$] (under a binomial model with a uniform prior distribution), which is remarkably close to DHM's prediction of (approximately) 51\%.

Evidence in favor of the DHM model's same-side bias prediction against the absence of the same-side bias can be evaluated using an informed Bayesian binomial test with $k$ same-side outcomes out of $N$ trials, $k   \sim \text{Binomial}(\beta, N)$, assuming that the coin flips are independently and identically distributed across people and coins. We specified two competing hypotheses via the binomial success parameter $\beta$, where success denotes the coin landing on the same side it started from:
\begin{align}
    \label{eq:overal_same_side}
    \text{No same-side bias, } \mathcal{H}_{0}                   &: \beta = 0.5 \\ \nonumber
    \text{DHM same-side bias, } \mathcal{H}_{1} &: \beta \sim\text{Beta}(5100, 4900)_{[0.5, 1]}.
\end{align}
The highly informed $\text{Beta}(5100, 4900)_{[0.5, 1]}$ prior distribution is meant to adequately represent the DHM hypothesis of the same-side bias (\citealp{diaconis2007dynamical}; see Appendix~\ref{app:priors} for more details about the prior distribution settings and prior sensitivity analysis). The evidence is quantified by the Bayes factor \citep{etz2017haldane, jeffreys1935some, jeffreys1939theory, kass1995bayes}:
\begin{equation*}
    \text{BF}_{10} = \frac{p(\text{data} \mid  \mathcal{H}_{1})}{p(\text{data} \mid  \mathcal{H}_{0})},
\end{equation*}
which contrasts the competing hypotheses in terms of their predictive performance for the observed data. The Bayes factor hypothesis test indicates extreme evidence in favor of the same-side bias predicted by the DHM model, $\text{BF}_{\text{same-side bias}} = 1.76\times10^{17}$. 

\begin{table}[h]
\caption{By-person summary of the probability of a same side landing.}
 \small
 \label{tab:people}
\centering
\begin{tabular}{lrrrrl}
  Person & Same side & Flips & Coins & Proportion [95\% CI] & Joined \\ 
  \hline
 XiaoyiL        &  780 &  1600 &    2 & 0.487 [0.463, 0.512] & Marathon-MSc \\ 
 JoyceP         & 1126 &  2300 &    3 & 0.490 [0.469, 0.510] & Marathon-MSc \\ 
 AndreeaZ       & 2204 &  4477 &    4 & 0.492 [0.478, 0.507] & Marathon-MSc \\ 
 KaleemU        & 7056 & 14324 &    8 & 0.493 [0.484, 0.501] & Bc Thesis \\ 
 FelipeFV       & 4957 & 10015 &    3 & 0.495 [0.485, 0.505] & Internet \\ 
 ArneJ          & 1937 &  3900 &    4 & 0.497 [0.481, 0.512] & Marathon-MSc \\ 
 AmirS          & 7458 & 15012 &    6 & 0.497 [0.489, 0.505] & Bc Thesis \\ 
 ChrisGI        & 4971 & 10005 &    5 & 0.497 [0.487, 0.507] & Marathon-Manheim \\ 
 FrederikA      & 5219 & 10500 &    5 & 0.497 [0.487, 0.507] & Internet \\ 
 FranziskaN     & 5368 & 10757 &    3 & 0.499 [0.490, 0.508] & Internet \\ 
 JasonN         & 3352 &  6700 &    7 & 0.500 [0.488, 0.512] & Marathon-PhD \\ 
 RietvanB       & 1801 &  3600 &    4 & 0.500 [0.484, 0.517] & Marathon-PhD \\ 
 PierreG        & 7506 & 15000 &    9 & 0.500 [0.492, 0.508] & Bc Thesis \\ 
 KarolineH      & 2761 &  5500 &    5 & 0.502 [0.489, 0.515] & Marathon-PhD \\ 
 SjoerdT        & 2510 &  5000 &    5 & 0.502 [0.488, 0.516] & Marathon-MSc \\ 
 SaraS          & 5022 & 10000 &    3 & 0.502 [0.492, 0.512] & Marathon-Manheim \\ 
 HenrikG        & 8649 & 17182 &    8 & 0.503 [0.496, 0.511] & Marathon     \\ 
 IrmaT          &  353 &   701 &    1 & 0.504 [0.467, 0.540] & Bc Thesis    \\ 
 KatharinaK     & 2220 &  4400 &    5 & 0.504 [0.490, 0.519] & Marathon-PhD \\ 
 JillR          & 3261 &  6463 &    2 & 0.505 [0.492, 0.517] & Marathon     \\ 
 FrantisekB     & 10148& 20100 &   11 & 0.505 [0.498, 0.512] & Marathon     \\ 
 IngeborgR      & 4340 &  8596 &    1 & 0.505 [0.494, 0.515] & Marathon     \\ 
 VincentO       & 2475 &  4900 &    5 & 0.505 [0.491, 0.519] & Marathon-MSc \\ 
 EricJW         & 2071 &  4100 &    5 & 0.505 [0.490, 0.520] & Marathon-MSc \\ 
 MalteZ         & 5559 & 11000 &    7 & 0.505 [0.496, 0.515] & Marathon-Manheim \\ 
 TheresaL       & 1769 &  3500 &    4 & 0.505 [0.489, 0.522] & Marathon-MSc \\ 
 DavidV         & 7586 & 14999 &    5 & 0.506 [0.498, 0.514] & Bc Thesis    \\ 
 AntonZ         & 5069 & 10004 &    2 & 0.507 [0.497, 0.516] & Marathon-Manheim \\ 
 MagdaM         & 2510 &  4944 &    6 & 0.508 [0.494, 0.522] & Marathon-MSc \\ 
 ThomasB        & 2540 &  5000 &    5 & 0.508 [0.494, 0.522] & Marathon-PhD \\ 
 JonasP         & 5080 &  9996 &    5 & 0.508 [0.498, 0.518] & Marathon     \\ 
 BohanF         & 1118 &  2200 &    3 & 0.508 [0.487, 0.529] & Marathon-MSc \\ 
 HannahA        & 1525 &  3000 &    4 & 0.508 [0.490, 0.526] & Marathon-MSc \\ 
 AdrianK        & 1749 &  3400 &    3 & 0.514 [0.498, 0.531] & Marathon-MSc \\ 
 AaronL         & 3815 &  7400 &    5 & 0.515 [0.504, 0.527] & Marathon-MSc \\ 
 KoenD          & 3309 &  6400 &    7 & 0.517 [0.505, 0.529] & Marathon-PhD \\ 
 MichelleD      & 2224 &  4300 &    5 & 0.517 [0.502, 0.532] & Marathon-PhD \\ 
 RoyMM          & 2020 &  3900 &    4 & 0.518 [0.502, 0.534] & Marathon-MSc \\ 
 TingP          & 1658 &  3200 &    4 & 0.518 [0.501, 0.535] & Marathon-MSc \\ 
 MaraB          & 1426 &  2750 &    3 & 0.518 [0.500, 0.537] & Marathon-MSc \\ 
 AdamF          & 4335 &  8328 &    2 & 0.520 [0.510, 0.531] & Marathon     \\ 
 AlexandraS     & 9080 & 17434 &    8 & 0.521 [0.513, 0.528] & Marathon     \\ 
 MadlenH        & 3705 &  7098 &    1 & 0.522 [0.510, 0.534] & Marathon     \\ 
 DavidKL        & 7895 & 15000 &    1 & 0.526 [0.518, 0.534] & Bc Thesis    \\ 
 XiaochangZ     & 1869 &  3481 &    4 & 0.537 [0.520, 0.553] & Marathon-MSc \\ 
 FranziskaA     & 2055 &  3800 &    4 & 0.541 [0.525, 0.557] & Marathon-MSc \\ 
 JanY           &  956 &  1691 &    2 & 0.565 [0.542, 0.589] & Marathon-MSc \\ 
 TianqiP        & 1682 &  2800 &    3 & 0.601 [0.582, 0.619] & Marathon-MSc \\ 
  \hline
  \textbf{Combined} & 178079 & 350757 & 44 & $0.508$ [$0.506$, $0.509$] &  \\
\end{tabular}
  \\
  \textit{Note.} `Proportion' refers to the observed proportion of coin flips that landed on the same side with a 95\% central credible interval under uniform prior distributions (virtually identical to a frequentist confidence interval).
\end{table}

\subsection*{Heads-tails bias}

Before proceeding to an analysis of heterogeneity in the same-side bias, we perform a similar analysis for the presence vs. absence of the heads-tails bias, i.e., $k \sim \text{Binomial}(\alpha, N)$ where $k$ corresponds to the number of heads out $N$ tosses with probability of landing on heads $\alpha$. Specifically, we obtained $175{,}421$ heads out of $350{,}757$ tosses, $\text{Pr}(\text{heads}) = 0.5001$, 95\% CI [$0.4985$, $0.5018$]. Replacing the DHM same-side bias prior distribution with a highly informed $\text{Beta}(5000, 5000)$ prior distribution (representing the background knowledge that in case a heads or tails bias exists, the effect is likely to be very small) yields moderate evidence against the presence of a heads-tails bias, $\text{BF}_{\text{heads-tails bias}} = 0.168$ (see Table~\ref{tab:coins} for by-coin summary). 

\begin{table}[h]
\caption{By-coin summary of the probability of heads.}
\label{tab:coins}
 \small
\centering
\begin{tabular}{lrrrr}
  Coin & Heads & Flips & People & Proportion [95\% CI]  \\ 
  \hline
  0.25CAD       &   48 &   100 &    1 & 0.480 [0.379, 0.582] \\ 
  20DEM-silver  &  484 &  1000 &    1 & 0.484 [0.453, 0.515] \\ 
  5CZK          & 1222 &  2500 &    2 & 0.489 [0.469, 0.509] \\ 
  0.05NZD       &  984 &  2011 &    1 & 0.489 [0.467, 0.511] \\ 
  0.10EUR       & 4515 &  9165 &    6 & 0.493 [0.482, 0.503] \\ 
  1DEM          & 2464 &  5000 &    5 & 0.493 [0.479, 0.507] \\ 
  50CZK         & 3207 &  6500 &    7 & 0.493 [0.481, 0.506] \\ 
  2HRK          & 4258 &  8596 &    1 & 0.495 [0.485, 0.506] \\ 
  1MXN          & 4180 &  8434 &    1 & 0.496 [0.485, 0.506] \\ 
  1SGD          & 7655 & 15400 &    2 & 0.497 [0.489, 0.505] \\ 
  5ZAR          & 3645 &  7325 &    1 & 0.498 [0.486, 0.509] \\ 
  2EUR          & 24276& 48772 &   28 & 0.498 [0.493, 0.502] \\ 
  0.01GBP       &  498 &  1000 &    1 & 0.498 [0.467, 0.529] \\ 
  0.50EUR       & 28617& 57445 &   32 & 0.498 [0.494, 0.502] \\ 
  0.20EUR       & 15665& 31373 &   20 & 0.499 [0.494, 0.505] \\ 
  0.25BRL       & 1998 &  4000 &    2 & 0.499 [0.484, 0.515] \\ 
  0.10RON       & 1000 &  2001 &    1 & 0.500 [0.478, 0.522] \\ 
  1CHF          & 2249 &  4500 &    4 & 0.500 [0.485, 0.514] \\ 
  1EUR          & 18920& 37829 &   25 & 0.500 [0.495, 0.505] \\ 
  0.20GEL       & 4501 &  8998 &    5 & 0.500 [0.490, 0.511] \\ 
  1CAD          & 5604 & 11200 &   11 & 0.500 [0.491, 0.510] \\ 
  2CAD          & 1502 &  3000 &    3 & 0.501 [0.483, 0.519] \\ 
  2MAD          & 1503 &  3000 &    1 & 0.501 [0.483, 0.519] \\ 
  100JPY        &  752 &  1500 &    1 & 0.501 [0.476, 0.527] \\ 
  2CHF          & 2259 &  4503 &    2 & 0.502 [0.487, 0.516] \\ 
  5MAD          & 1007 &  2001 &    1 & 0.503 [0.481, 0.525] \\ 
  0.20GBP       & 1516 &  3005 &    2 & 0.504 [0.486, 0.523] \\ 
  1CNY          &  757 &  1500 &    1 & 0.505 [0.479, 0.530] \\ 
  1CZK          &  505 &  1000 &    1 & 0.505 [0.474, 0.536] \\ 
  2ILS          &  506 &  1000 &    1 & 0.506 [0.475, 0.537] \\ 
  5JPY          & 1772 &  3500 &    2 & 0.506 [0.490, 0.523] \\ 
  5SEK          & 8052 & 15902 &    7 & 0.506 [0.499, 0.514] \\ 
  0.25USD       & 2180 &  4300 &    4 & 0.507 [0.492, 0.522] \\ 
  1MAD          & 1014 &  2000 &    1 & 0.507 [0.485, 0.529] \\ 
  0.50RON       & 1442 &  2844 &    3 & 0.507 [0.488, 0.526] \\ 
  0.05EUR       & 3821 &  7514 &    6 & 0.509 [0.497, 0.520] \\ 
  0.50GBP       &  765 &  1504 &    1 & 0.509 [0.483, 0.534] \\ 
  2BDT          & 2038 &  4003 &    2 & 0.509 [0.494, 0.525] \\ 
  10CZK         & 4572 &  8905 &    7 & 0.513 [0.503, 0.524] \\ 
  0.20CHF       &  518 &  1000 &    1 & 0.518 [0.487, 0.549] \\ 
  0.50SGD       & 1449 &  2781 &    3 & 0.521 [0.502, 0.540] \\ 
  0.02EUR       &  158 &   300 &    1 & 0.527 [0.468, 0.584] \\ 
  1GBP          &  791 &  1500 &    2 & 0.527 [0.502, 0.553] \\ 
  2INR          &  552 &  1046 &    1 & 0.528 [0.497, 0.558] \\ 
   \hline
  \textbf{Combined} & 175421 & 350757 & 48 & $0.500$ [$0.498$, $0.502$] \\ 
\end{tabular}
  \\
  \textit{Note.} `Proportion' refers to the observed proportion of coin flips that landed heads with a 95\% central credible interval under uniform prior distributions (virtually identical to a frequentist confidence interval).
\end{table}

\subsection*{Heterogeneity: Between-people and between-coins}

A closer examination of the individual results, however, suggests that this overall result needs to be qualified in the sense that the size of the same-side bias varies across individuals (cf. Table~\ref{tab:people}). We extend the binomial model into a hierarchical model with by-person $k = 1, \dots, K$ and by-coin $j = 1, \dots, J$ specific deviations $\gamma_{\beta_k}$ and $\gamma_{\alpha_j}$ from the overall same-side and the heads-tails bias (i.e., $\text{logit}(\beta_\mu)$ and $\text{logit}(\alpha_\mu)$, respectively\footnote{While $\beta_\mu$ and $\alpha_\mu$ correspond to the overall probability of the same-side and heads respectively, the logit() function transforms the overall probability into the bias (e.g., equal probability of heads and tails, $\alpha_\mu = 0.5$, corresponds to no heads-tails bias, $\text{logit}(\alpha_\mu) = 0$).}):
\begin{align}
    \label{eq:model}
    \gamma_{\alpha_j}       &\sim \text{Normal}(0, \sigma_\alpha^2)              \\  \nonumber
    \gamma_{\beta_k}        &\sim \text{Normal}(0, \sigma_\beta^2)               \\  \nonumber
    \text{logit}(\alpha_j)  &=    \underbrace{\text{logit}(\alpha_\mu)}_{\text{Overall same-side bias}} + \underbrace{\gamma_{\alpha_j}}_{\substack{\text{Person-specific deviation}\\\text{from same-side bias}}}   \\  \nonumber
    \text{logit}(\beta_k)   &=    \underbrace{\text{logit}(\beta_\mu)}_{\text{Overall heads-tails bias}} + \underbrace{\gamma_{\beta_k}}_{\substack{\text{Coin-specific deviation}\\\text{from heads-tails bias}}}     \\  \nonumber
    \mu_{ijk}      &=  
    \begin{cases}
        \text{logit}(\alpha_j) + \text{logit}(\beta_k)  \,\,\,\,   & \text{y}_{s=0,ijk} = 1 \,\text{(Starting heads)} \\ 
        \underbrace{\text{logit}(\alpha_j)}_{\substack{\text{Heads bias}}} - \underbrace{\text{logit}(\beta_k)}_{\substack{\text{Same-side}\\\text{bias}}}  \,\,\,\,   & \text{y}_{s=0,ijk} = 0 \,\text{(Starting tails)} \\ 
    \end{cases}\\
    \text{y}_{s=1,ijk}      &\sim \text{Bernoulli}(\text{logit}^{-1}(\mu_{ijk})), \nonumber
\end{align}
where $\text{y}_{s=1,ijk} = 1$ corresponds to the $i^\text{th}$ flip of the $k^\text{th}$ person with the $j^\text{th}$ coin landing heads.\footnote{Since the DHM model specifies that the same-side bias results from the wobble introduced by the coin tosser (i.e., the person), the model does not feature coin-specific random-effects in the same-side bias (cf. Appendix C for results suggesting no between-coin heterogeneity in same-side bias). Similarly, the model does not feature person-specific random-effects in heads-tails bias as we would not expect that some people are more likely to flip heads than tails or visa versa.} We first use the model with slightly informative prior distributions tailored for parameter estimation, and then proceed with more informed prior distributions tailored for hypothesis testing and Bayesian model-averaging to test simultaneously for the presence vs. absence of the same-side bias, heads-tails bias, and between-people and between-coin heterogeneities in the respective biases  \citep[cf.,][]{jeffreys1939theory, jeffreys1961theory}. See Appendix~\ref{app:priors} for more details.

\subsubsection*{Parameter estimation}

\begin{figure}[h]
    \centering
    \includegraphics[width=\textwidth]{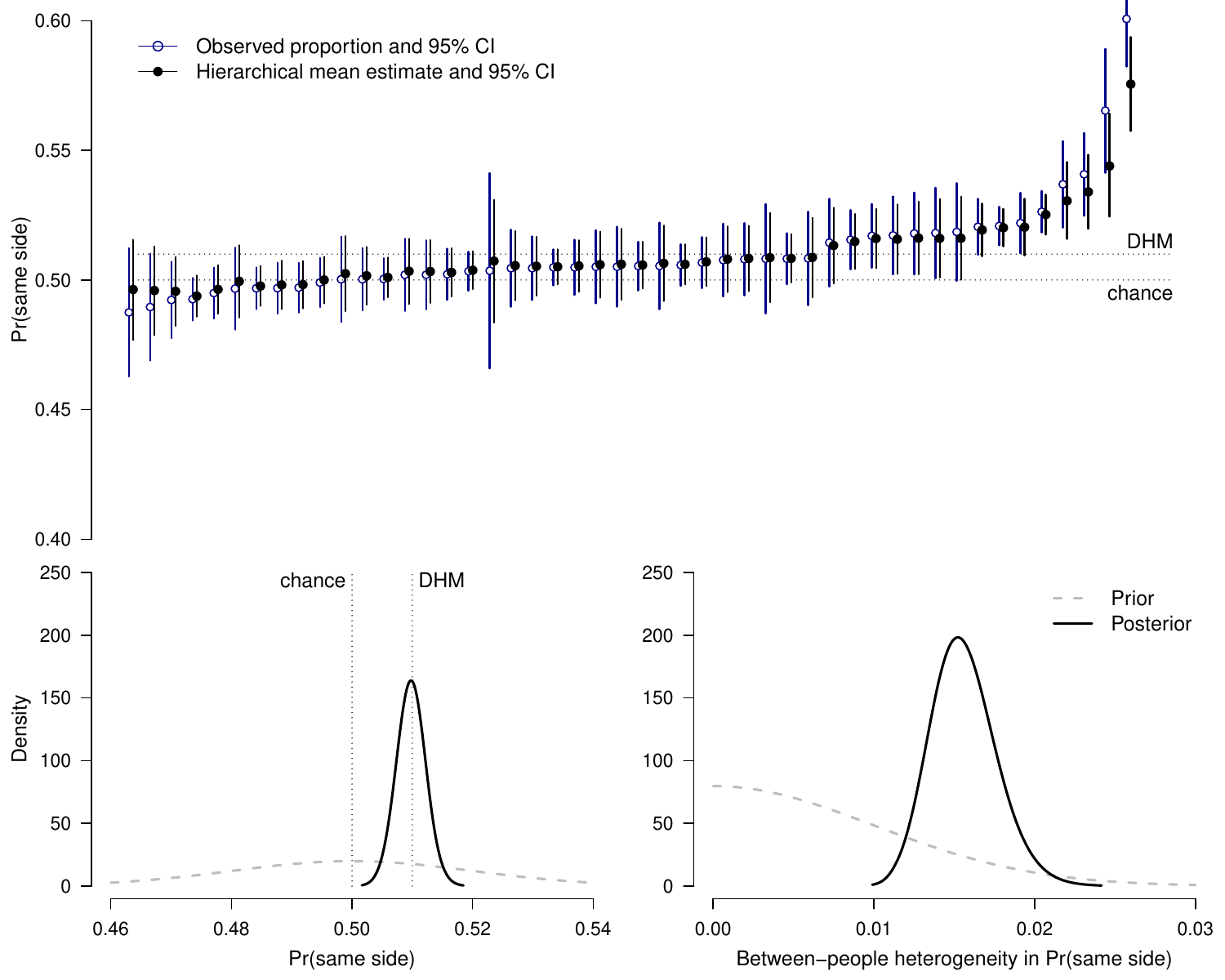}
    \caption{Coins have a tendency to land on the same side they started, confirming the predictions from the Diaconis, Holmes, and Montgomery (DHM) model of coin flipping. Top panel: posterior estimates of the probability of same side separately for each person, as obtained from the hierarchical Bayesian model with weakly informative, estimation-tailored prior distributions described in the methods section; Bottom-left panel: prior and posterior distributions for the overall probability of same side; Bottom-right panel: prior and posterior distributions for the between-people heterogeneity in the probability of the same side.}
    \label{fig:people}
\end{figure}

The observed proportion (and 95\% CI) of same-side outcomes for each individual participant is shown as the blue dots (and bars) in Figure~\ref{fig:people}. It is clear that participants notably differ in the degree of the same-side bias. The black dots correspond to the estimated probability of the same-side outcome from the hierarchical model with slightly informative prior distributions that incorporates both between-people heterogeneity in same-side bias and between-coin heterogeneity in heads-tails bias. Only a small proportion of participants' point estimates fall below the chance line (i.e., no bias), with the majority of point estimates somewhere in between the chance and the DHM model prediction, and some exceeding the DHM prediction. The left bottom panel of Figure~\ref{fig:people} displays the posterior distribution of the probability of a same side outcome, $\text{Pr}(\text{same side}) = 0.5098$, 95\% CI [$0.5050$, $0.5147$]. This credible interval is slightly wider than the one from the simple binomial analysis; this is due to the substantial between-people heterogeneity in the probability of the coin landing on the same side, $\text{sd}_\text{people}(\text{Pr}(\text{same side})) = 0.0156$, 95\% CI [$0.0119$, $0.0200$] (right bottom panels of Figure~\ref{fig:people}).

\begin{figure}[h]
    \centering
    \includegraphics[width=1\textwidth]{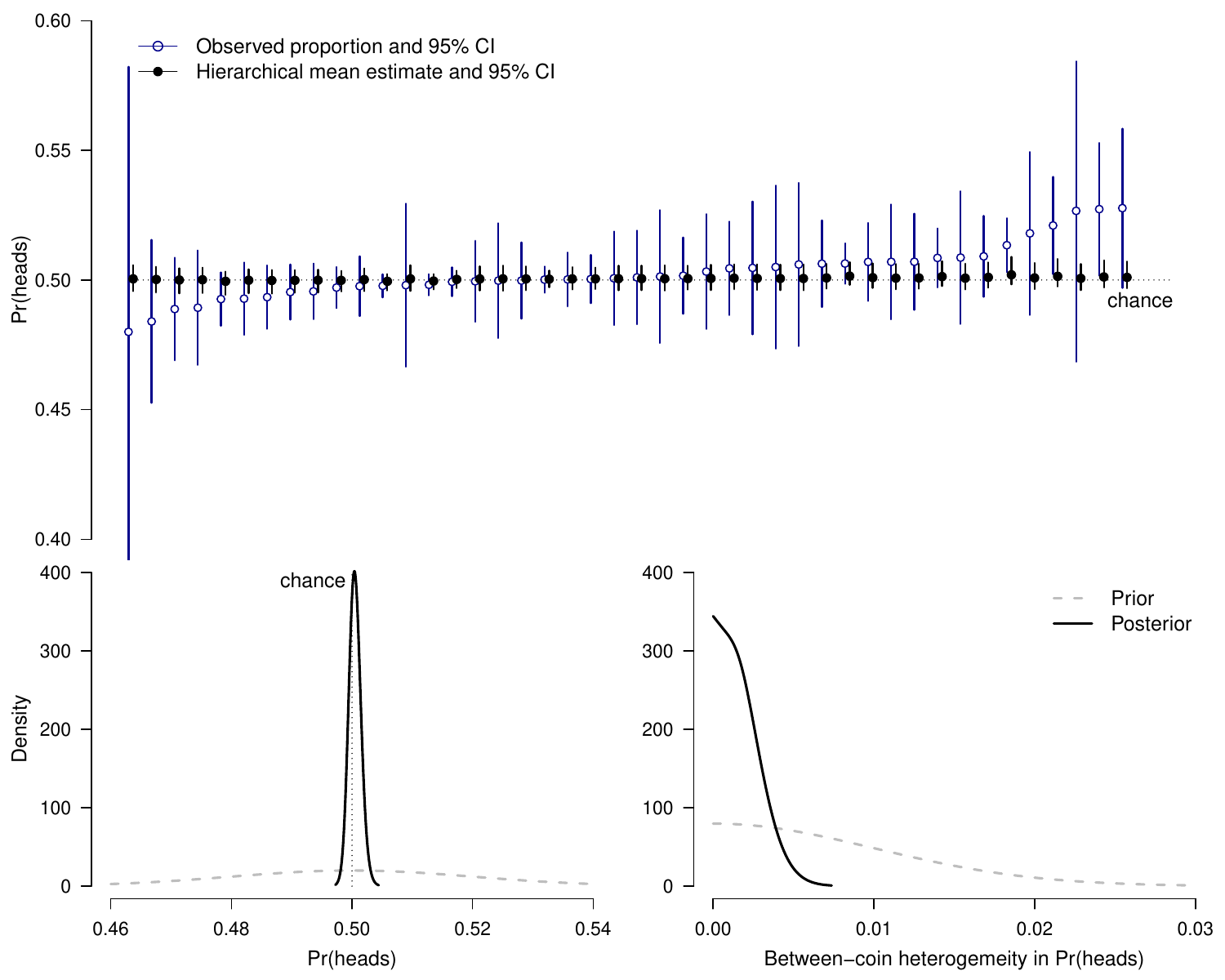}
    \caption{Coins have a tendency to land on heads and tails with equal probability, supporting the predictions from the Standard model of coin flipping. Top panel: posterior estimates of the probability of heads separately for each coin, as obtained from the hierarchical Bayesian model with weakly informative, estimation-tailored prior distributions described in the methods section; Bottom-left panel: prior and posterior distributions for the overall probability of heads; Bottom-right panel: prior and posterior distributions for the between-coin heterogeneity in the probability of heads.}
    \label{fig:coins}
\end{figure}

Inspection of the coin-specific parameters (Figure~\ref{fig:coins}) suggests the lack of a heads-tails bias, $\text{Pr}(\text{heads}) = 0.5005$, 95\% CI [$0.4986$, $0.5026$], with virtually no between-coin heterogeneity, $\text{sd}_\text{coins}(\text{Pr}(\text{heads})) = 0.0018$, 95\% CI [$0.0001$, $0.0047$], as all observed proportions and estimated probabilities of landing on heads cluster around the chance line.

\subsubsection*{Hypothesis testing}
Evidence in favor of the presence vs. absence of the same-side bias, heads-tails bias, and between-people and between-coin heterogeneities can be evaluated simultaneously using Bayesian model-averaging \citep{raftery1995accounting, hoeting1999bayesian, hinne2019conceptual} and \emph{inclusion Bayes factors}, a generalization of Bayes factors based on the change from prior to posterior odds \citep{hinne2019conceptual}:
\begin{equation}
    \label{eq:inclusion-BF}
\underbrace{ \text{BF}_{\text{AB}}}_{ \substack{\text{Inclusion Bayes factor}\\{\text{for A vs. B}}} } =  
    \;\;\;  \underbrace{ \frac{ \sum_{a \in A} p(\mathcal{M}_{a} \mid \text{data}) }
    { \sum_{b \in B} p(\mathcal{M}_{b} \mid \text{data}) }}_{ \substack{\text{Posterior inclusion odds}\\{\text{for A vs. B}}}} \;\;\; \Bigg/ \underbrace{ \frac{ \sum_{a \in A} p(\mathcal{M}_{a}) }
    { \sum_{b \in B} p(\mathcal{M}_{b}) }}_{\substack{\text{Prior inclusion odds}\\{\text{for A vs. B}}}},
\end{equation}
where $A$ contains a set of models where a given hypothesis holds and $B$ contains the complement. The overarching set of all models ($A$ and $B$) can be specified as an orthogonal combination of the different possible hypotheses, i.e., $2^4 = 16$ models, where each is assigned the usual equal prior model probability of $1/16$ (Appendix~\ref{app:priors} lists all specified models in detail).

The inclusion Bayes factor indicated compelling evidence for the presence of an overall same-side bias, $\text{BF}_{\text{same-side bias}} = 2359$. When the hypothesis of a same-side bias has a prior probability of $0.50$, a Bayes factor of about $2359$ results in a posterior probability of $0.9996$. In addition, consistent with the visual impression from Figure~\ref{fig:people}, the inclusion Bayes factor reveals overwhelming evidence for the presence of between-people heterogeneity in same-side bias, $\text{BF}_{\text{people heterogeneity}} = 3.10 \times 10^{24}$. The evidence against the presence of heads-tails bias remains practically unchanged, $\text{BF}_{\text{heads-tails bias}} = 0.182$. The model indicates moderate evidence against the presence of between-coin heterogeneity in heads-tails bias, $\text{BF}_{\text{coin heterogeneity}} = 0.178$ (see Appendix C for concordant frequentist result and Table 1 in the Online Supplements for an overview of the individual models).

\subsection*{Practice effects}
\begin{figure}[h]
    \centering
    \includegraphics[width=1\textwidth]{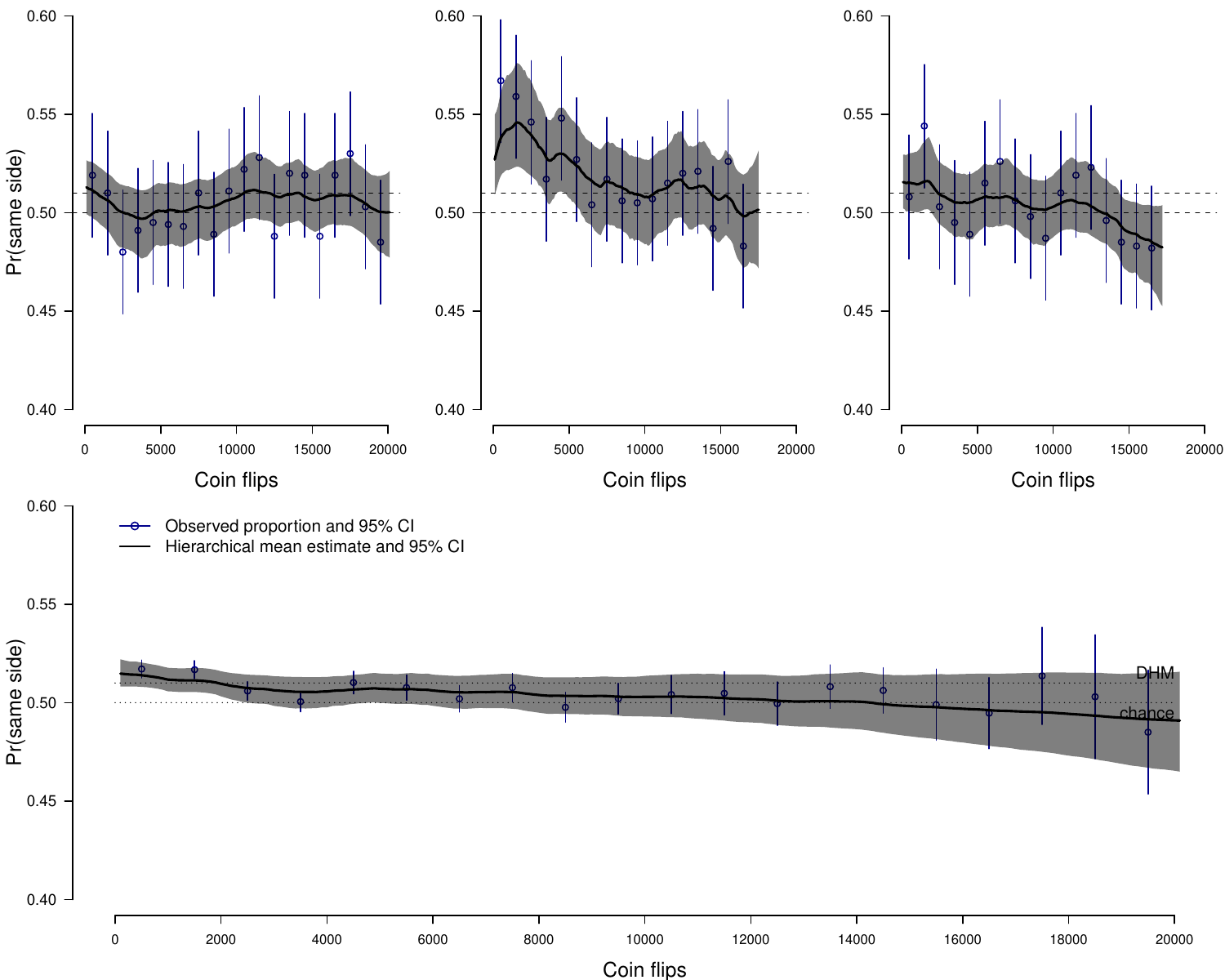}
    \caption{The probability of coins landing on the same side changes with practice. Top panel: posterior estimates of the probability of the same side separately for three participants with the largest number of coin flips across the number of conducted coin flips (see Online Supplements for all trajectories); Bottom panel: posterior estimates of the probability of the same side across all participants showing on average decreasing pattern of the probability of the same side across the number of conducted coin flips. The posterior estimates are obtained from the hierarchical Bayesian state-space model with weakly informative, estimation-tailored prior distributions described in the methods section overlaying the observed proportions binned by 1000 coin flips.}
    \label{fig:time}
\end{figure}

Following up on a suggestion by a reviewer, a more extensive inspection of the data suggested that the degree of same-side bias changes with the number of coin flips performed (i.e., practice) by each person. Moreover, the pattern and degree of change varied across participants. The top panel of Figure~\ref{fig:time} shows the proportion of same-side outcomes aggregated by 1000 coin flips as black dots (and bars) for three participants with the largest number of coin flips. Each panel displays a qualitatively different pattern of change; the left panel shows a wavy pattern of same-side bias with practice, the middle panel shows a pattern in which the same-side bias first decreases and then remains stationary, and the right panel shows a pattern in which the same-side bias is first stationary and then decreases. While not all individual patterns fit into one of these three ``pattern types'', most participants showed an initial decrease in same-side bias with practice (see Online Supplements for visualization of all participant trajectories). The bottom panel of Figure~\ref{fig:time} shows the proportion of same-side outcomes aggregated by 1000 coin flips as black dots (and bars) combined across all participants. The overall pattern of same-side bias seems to be monotonically decreasing with practice. Note that the uncertainty increases at the higher numbers of coin flips, because only a few participants flipped coins more than 15,000 times.

The black trend lines and credible bands in Figure~~\ref{fig:time} correspond to the estimated probability of a same-side outcome from a hierarchical state-space model that monitors the change in the same-side bias with practice. State-space models are flexible non-parametric models often used for monitoring non-monotonic and non-linear changes over time \citep{fruhwirth2006finite}. Here, the hierarchical state-space model extends the hierarchical binomial model (Eq.~\ref{eq:model}) by modeling the current state of the same-side bias of the $i^\text{th}$ flip of the $k^\text{th}$ person as a person-and-flip-specific change $\delta_{ki}$ from their previous state of the same-side bias, $\beta_{k(i-1)})$, and a constant person-specific drift parameter $\theta_k$,

\begin{align}
    \label{eq:modelt}
    \gamma_{\beta_k}         &\sim \text{Normal}(0, \sigma_\beta^2)               \\  \nonumber
    \gamma_{\delta_{ki}}     &\sim \text{Normal}(0, \sigma_\delta^2)              \\  \nonumber
    \gamma_{\theta_k}        &\sim \text{Normal}(0, \sigma_\theta^2)              \\  \nonumber
    \delta_{ki}              &=   \underbrace{\delta_{\mu i}}_{\substack{\text{Overall state change}\\\text{at flip $i$}}} + \underbrace{\gamma_{\delta_{ki}}}_{\substack{\text{Person-specific deviation}\\\text{from state change at flip $i$}}}          \\  \nonumber
    \theta_k                 &=   \underbrace{\theta_\mu}_{\text{Overall drift}} + \underbrace{\gamma_{\theta_k}}_{\substack{\text{Person-specific deviation}\\\text{from drift}}}                \\  \nonumber
    \underbrace{\text{logit}(\beta_{ki})}_{\substack{\text{Current state} \\ \text{of same-side bias}}} &= 
    \begin{cases}
             \underbrace{\text{logit}(\beta_{\mu 1})}_{\substack{\text{Overall initial} \\ \text{same-side bias}}} \,\,\,\,\, + \underbrace{\gamma_{\beta_{k 1}}}_{\substack{\text{Initial person-specific deviation} \\ \text{from same-side bias}}}   & i = 1 \, \text{(Initial state)}\\ 
             \underbrace{\text{logit}(\beta_{k(i-1)})}_{\substack{\text{Previous state} \\ \text{of same-side bias}}} \,\,\,\,\, + \underbrace{\delta_{ki}}_{\substack{\text{State change}\\\text{at flip $i$}}} + \underbrace{\theta_k}_{\text{Drift}} & i > 1 \,  \text{(Subsequent states)} \\ 
     \end{cases}\\
    \mu_{ijk}      &=  
    \begin{cases}
        \,\,\,\,\,\text{logit}(\beta_{ki})  \,\,\,\,   & \text{y}_{s=0,ijk} = 1 \,\text{(Starting heads)} \\ 
        - \text{logit}(\beta_{ki})  \,\,\,\,   & \text{y}_{s=0,ijk} = 0 \,\text{(Starting tails)} \\ 
    \end{cases}\\
    \text{y}_{s=1,ijk}      &\sim \text{Bernoulli}(\text{logit}^{-1}(\mu_{ijk})), \nonumber
\end{align}
where $\delta_{\mu i}$ corresponds to the overall same-side bias at $i^\text{th}$ flip and $\theta_\mu$ corresponds to the overall drift in the same-side bias. In contrast to the earlier model, Eq.~\eqref{eq:model}, we omit the overall heads-tails bias and the between-coin heterogeneity in the heads-tails bias to simplify the model estimation process in light of the evidence against the presence of the heads-tails bias and its heterogeneity. We only use the model with slightly informative prior distributions tailored for parameter estimations, see Appendix~\ref{app:priors} for more details.

The pooled estimate of the initial same-side bias based on the state-space model equals $\text{Pr}(\text{same side}) = 0.5148$, 95\% CI [$0.5082$, $0.5220$] and is larger than the overall same-side bias estimate from the hierarchical model. However, the initial heterogeneity in the probability of the coin landing on the same side (based on the state-space model; $\text{sd}_\text{people}(\text{Pr}(\text{same side})) = 0.0082$, 95\% CI [$0.0008$, $0.0162$]) is almost twice as low as the overall heterogeneity in the probability of coins landing on the same side (based on the hierarchical model; $\text{sd}_\text{people}(\text{Pr}(\text{same side})) = 0.0156$, 95\% CI [$0.0119$, $0.0200$]). This result suggests the large between-people heterogeneity observed in the hierarchical model was confounded to some degree by practice effects. 

\subsection*{Outlier exclusion}

We repeated the statistical analyses after excluding four potential outliers with same-side sample proportions larger than 53\% (i.e., the four largest and right-most estimates in the top panel of Figure~\ref{fig:people}). In general, the exclusion did not qualitatively affect our conclusions, although --as may be expected-- the same-size bias decreased in size, and the between-people heterogeneity became less pronounced. Each of the four participants with the same-side sample proportion larger than 53\% contributed fewer than $4{,}000$ coin flips, which is consistent with the possibility that the relatively high proportions may be due in part to the fact that the same-side bias is largest at the beginning (before the practice effects occur). Additional robustness checks reported in Tables~3, 4, and 5 in the Online Supplements demonstrate that the qualitative conclusions do not change when excluding 1 to 5 participants with the lowest and the largest proportion of same-side outcomes.

After excluding the four potential outliers, the data feature $171{,}517$ same-side landings from $338{,}985$ tosses, $\text{Pr}(\text{same side}) = 0.5060$, 95\% CI [$0.5043$, $0.5077$]. The evidence in favor of the DHM hypothesis using the Bayesian informed binomial hypothesis test decreased notably but remains extreme, $\text{BF}_{\text{same-side bias}} = 1.28\times10^{8}$. The proportion of heads remained practically identical, $169{,}635$ heads out of $338{,}985$ tosses, $\text{Pr}(\text{heads}) = 0.5004$, 95\% CI [$0.4987$, $0.5021$], as did the moderate evidence against heads-tails bias $\text{BF}_{\text{heads-tails bias}} = 0.190$.

The exclusion of potential outliers similarly affects the inference from the hierarchical model: although the same-side bias decreases and the associated heterogeneity is reduced (i.e., $\text{Pr}(\text{same side}) = 0.5060$, 95\% CI [$0.5031$, $0.5089$] and $\text{sd}_\text{people}(\text{Pr}(\text{same side})) = 0.0072$, 95\% CI [$0.0050$, $0.0099$]), the evidence for the presence of the same-side bias and the associated between-people heterogeneity remains extreme (i.e., $\text{BF}_{\text{same-side bias}} = 787$ and $\text{BF}_{\text{people heterogeneity}} = 2.87 \times 10^{7}$). For the probability of heads, the inference remains practically unchanged, both with respect to the size of the effect and associated heterogeneity (i.e., $\text{Pr}(\text{heads}) = 0.5008$, 95\% CI [$0.4988$, $0.5030$] and $\text{sd}_\text{coins}(\text{Pr}(\text{heads})) = 0.0020$, 95\% CI [$0.0001$, $0.0050$]), and with respect to the evidence for the presence of the heads-tails bias and the associated between-coin heterogeneity (i.e., $\text{BF}_{\text{heads-tails bias}} = 0.213$ and $\text{BF}_{\text{coin heterogeneity}} = 0.221$).

The exclusion of potential outliers affects the practice effects model in a similar manner. The initial same-side bias slightly decreased, $\text{Pr}(\text{toss-order same side}) = 0.5104$, 95\% CI [$0.5049$, $0.5167$], and the associated between-people heterogeneity was reduced, $\text{sd}_\text{people}(\text{Pr}(\text{toss-order same side})) = 0.0036$, 95\% CI [$0.0003$, $0.0099$].

\subsection*{Including the Results from Larwood and Ku}

We repeated the statistical analyses after including the $40{,}000$ coin flips performed by Larwood and Ku \cite{noauthor_40000_nodate}. In general, the inclusion did not qualitatively affect our conclusion, although the evidence for the same-side bias slightly increased.

After inclusion of the $40{,}000$ coin flips, the data feature ${198,324}$ same-side landings from ${390,757}$ tosses, $\text{Pr}(\text{same side}) = 0.5075$, 95\% CI [$0.5060$, $0.5091$]. The evidence in favor of the DHM hypothesis using the Bayesian informed binomial hypothesis test slightly increases, $\text{BF}_{\text{same-side bias}} = 2.77\times10^{18}$. The proportion of heads remained practically identical, ${195,638}$ heads out of $390,757$ tosses, $\text{Pr}(\text{heads}) = 0.5007$, 95\% CI [$0.4991$, $0.5022$], as did the moderate evidence against heads-tails bias $\text{BF}_{\text{heads-tails bias}} = 0.221$.

The inclusion of the $40{,}000$ coin flips affects the inference from the hierarchical model in a similar fashion: the degree of the same-side bias and the associated heterogeneity remains practically unchanged (i.e., $\text{Pr}(\text{same side}) = 0.5096$, 95\% CI [$0.5050$, $0.5142$] and $\text{sd}_\text{people}(\text{Pr}(\text{same side})) = 0.0152$, 95\% CI [$0.0116$, $0.0194$]), the evidence for the presence of the same-side bias increases, and the evidence for the associated between-people heterogeneity slightly decreases (i.e., $\text{BF}_{\text{same-side bias}} = 3295$ and $\text{BF}_{\text{people heterogeneity}} = 1.66 \times 10^{24}$). For the probability of heads, the inference also remains practically unchanged, both with respect to the size of the effect and associated heterogeneity (i.e., $\text{Pr}(\text{heads}) = 0.5005$, 95\% CI [$0.4986$, $0.5026$] and $\text{sd}_\text{coins}(\text{Pr}(\text{heads})) = 0.0017$, 95\% CI [$0.0001$, $0.0047$]), and with respect to the evidence for the presence of the heads-tails bias and the associated between-coin heterogeneity (i.e., $\text{BF}_{\text{heads-tails bias}} = 0.183$ and $\text{BF}_{\text{coin heterogeneity}} = 0.175$).

The inclusion of the $40{,}000$ coin flips similarly affects the practice effects model. The initial same-side bias slightly decreased, $\text{Pr}(\text{toss-order same side}) = 0.5132$, 95\% CI [$0.5074$, $0.5195$], and the associated between-people heterogeneity was reduced, $\text{sd}_\text{people}(\text{Pr}(\text{toss-order same side})) = 0.0071$, 95\% CI [$0.0006$, $0.0148$].

\section*{Discussion}

We collected $350{,}757$ coin flips and found strong empirical evidence for the counterintuitive and precise prediction from DHM model of human coin tossing: when people flip a coin, it tends to land on the same side as it started. However, this conclusion needs to be qualified by two important factors. Firstly, the data revealed a substantial degree of between-people variability in the same-side bias: as can be seen in Figure~\ref{fig:people}, some people appear to have little or no same-side bias, whereas others do display a same-side bias, albeit to a varying degree. This variability is arguably consistent with DHM model, in which the same-side bias originates from off-axis rotations (i.e., precession or wobbliness), which can reasonably be assumed to vary between people. Secondly, additional exploratory analyses suggest that the degree of the same-side bias decreases with the number of coin flips. A possible explanation of this decreasing bias is a coin-tossing practice effect---the more coins people flip, the closer they approach the `perfect' wobble-less flip. Our results suggest that around $10{,}000$ coin flips ($\approx 10$ hours of coin flipping) might be enough to virtually eliminate the same-side bias. Our results are robust to the type of coins used and to changes in analytic methodology -- for instance, different prior distributions as outlined in Appendix~B, different analytical approaches as outlined in Appendix~D and discussed by Viechtbauer \cite{viechtbauerVignette} and Pawel \cite{pawel2024bayes}, and different criteria for outlier exclusion as outlined in the Supplementary Materials. However, our experiment cannot rule out alternative explanations of the reduction in the same-side bias over time, such as exhaustion, waning attention, etc.

Our results are aligned with previous empirical data evaluating same-side bias in coin flipping; the Diaconis et al. \cite{diaconis2007dynamical} 27 high-speed camera flips that initially suggested the same-side bias of approximately $1$\%, and the subsequent $20{,}000$ coin flips by Janet Larwood and Priscilla Ku that resulted in a same-side bias of $1.2$\% and $0.1$\% , respectively (results that incidentally also highlight the between-people heterogeneity in the same-side bias; \citealp{noauthor_40000_nodate}). Furthermore, in a recent informal replication, McGaw \cite{mcgaw2024} collected 10 coin flips from 820 people each, resulting in 4171 same-side outcomes, for a same-side proportion of 0.5087.

Future work may attempt to verify whether `wobbly tossers' show a more pronounced same-side bias than `stable tossers' and examine the practice effects in more detail. Furthermore, the present practice effects suggest that the same-side bias is best studied through many people who each contribute only a few thousand coin flips rather than through few people who each contribute tens of thousands of coin flips; the same-side bias is much more pronounced at the beginning of the experiment, and hence there are diminishing returns when a single person contributes very many coin flips. However, the effort required to test the more detailed hypotheses appears to be excessive, as this would ideally involve detailed analyses of high-speed camera recordings for individual flips (cf. \citealp{diaconis2007dynamical}).  

In order to ensure the quality of the data, we videotaped and audited the data collection procedure (see the Appendix section for details). The audit did not reveal anything suspicious (i.e., all participants performed the coin flips they reported). While the video recordings were of insufficient quality to provide additional insights into the variability of the same-side bias, McGaw \cite{mcgaw2024} noticed that the coin flipping technique of participants with the extreme same-side bias was visibly worse than that of the remaining participants. Exclusion of those participants, however, did not significantly affect the results. There also remains a legitimate concern: at the time when people were flipping the coins they were aware of the main hypothesis under test. Therefore it cannot be excluded that some of the participants were able to manipulate the coin flip outcomes in order to produce the same-side bias. In light of the nature of the coin tossing process, the evidence from the video recordings, and the precise correspondence between the data and the predictions from DHM model, we deem this possibility unlikely, but future work is needed to disprove it conclusively (e.g., by concealing the aim of the study). 

Could future coin tossers use the same-side bias to their advantage? The magnitude of the observed bias can be illustrated using a betting scenario. If you bet a dollar on the outcome of a coin toss (i.e., paying 1 dollar to enter, and winning either 0 or 2 dollars depending on the outcome) and repeat the bet $1{,}000$ times, knowing the starting position of the coin toss would earn you $19$ dollars on average. This is more than the casino advantage for $6$ deck blackjack against an optimal-strategy player, where the casino would make $5$ dollars on a comparable bet, but less than the casino advantage for single-zero roulette, where the casino would make $27$ dollars on average \citep{hannum2003guide}. Moreover, this advantage amounts to a non-negligible 95,000\$ in expected value when considering the 5,000,000\$ double-or-eliminated coin flip in the ``Beast Games'' Amazon Prime show \citep{beastgames2025}. These considerations lead us to suggest that when coin flips are used for high-stakes decision-making, the starting position of the coin is best concealed.

\section*{Acknowledgments}

\subsection*{Funding}
The authors have no funding to declare, and conducted this research in their spare time.


\subsection*{Competing interests}
The authors have no competing interest to declare.

\subsection*{Data and materials availability}
All data and materials are available at \url{https://osf.io/pxu6r/}.

\subsection*{Ethical approval}
The research project was approved by the Ethics Review Board of the Faculty of Social and Behavioral Sciences, University of Amsterdam, The Netherlands (2022-PML-15687).

\subsection*{Author note}
During a revision we corrected two minor errors in the data set: (a) two coins were miscoded as different currency $\times$ denomination, and (b) one invalid coin flip with unknown starting position was not removed from the analysis as incorrectly treated as starting heads and landing tails. We fixed the coin coding issue, which reduced the number of unique coins from 46 to 44. We fixed the mistake in data processing and asked the participant to provide one additional valid coin flip to retain the same number of total coin flips. The additional flip started tails and landed tails, which increased the number of coins landing on the same side by one and decreased the number of coins landing on heads by one (we could not retrieve the original coin; consequently, the number of flips with 5 ZAR reduced by one and number of flips with 0.05 EUR increased by one). Our original findings were qualitatively unaffected by these errors.

\bibliographystyle{WileyNJD-AMA}
\bibliography{manuscript}

\section*{Appendix A: Prior Distributions}
\label{app:priors}

\subsection*{Same-side bias}
The informed prior distribution for the same-side bias, $\beta \sim \text{Beta}(5100, 4900)$, was constructed to take into account the relevant background knowledge and obey several restrictions. Specifically: (a) the prior distribution is centered at the same-side bias of approximately 0.01 (the DHM point-prediction of about 51\%); (b) the prior distribution allows for reasonable variability around the mode---its standard deviation of approximately 0.005 corresponds to a 95\% prior CI on the same-side bias that ranges from 0.00 to 0.02 (aligned with the majority of the estimates in \citealp{diaconis2007dynamical}; however, note that \citealp{diaconis2007dynamical} also recorded a couple of outliers with a bias higher than 0.02); and (c) the prior distribution is truncation at 0.50, restricting the hypothesis to predict only the kind of bias stipulated by DHM (i.e., a same-side bias; a probability below 0.50 would represent an opposite-side bias, which would be anomalous in the context of the DHM model). Note that this informed hypothesis is motivated by the precise prediction of the DHM hypothesis---and as such allows for a meaningful test.

\subsection*{Heads-tails bias}
In the absence of a theory-driven prediction for the size of the heads-tails bias, we used a prior distribution that is similar to the one we employed for the same-side bias, $\alpha \sim \text{Beta}(5000, 5000)$, since: (a) if any heads-tails bias were present it would probably be minuscule---which is encoded in the 95\% prior CI ranging from 0.49 (i.e., a 0.01 bias for tails) to 0.51 (i.e., a 0.01 bias for heads), and (b) we lack a reason to prefer either heads or tails, and hence the prior distribution is symmetric around $\alpha = 0.5$.

\subsection*{Heterogeneity: Between-people and between-coins}
\subsubsection*{Parameter estimation}
For both the probability of heads $\alpha_\mu$ and the probability of the same-side $\beta_\mu$, we specified $\text{Beta}(312, 312)$ prior distributions. The $\text{Beta}(312, 312)$ prior distributions are centered at 0.50 (i.e., no bias) with a standard deviation of 0.02; we chose those values as the 95\% prior CI contains biases up to 4 times larger than those predicted by the DHM theory which grants sufficient flexibility while ruling out implausibly large biases. For both the between-coin heterogeneity $\sigma_\alpha$ and between-people heterogeneity $\sigma_\beta$ in the biases, we specified $\text{Normal}_{+}(0, 0.04)$ prior distributions. The $\text{Normal}_{+}(0, 0.04)$ prior distribution results in approximately half-Normal prior distributions with a standard deviation of 0.01 when transformed to probability scale; we chose those values as we believed that the between-coin and between-people heterogeneity in the biases would most likely be lower or equal to the bias predicted by DHM. Later visualizations confirmed that the specified prior distributions were not overly restrictive.

\subsubsection*{Hypothesis testing}
In order to test the alternative hypotheses of the same-side bias $\mathcal{H}_{\beta, 1}$ and the heads-tails bias $\mathcal{H}_{\alpha, 1}$ we specified beta prior distributions that were identical to those used in the non-hierarchical models. For the between-people and between-coin heterogeneity in the same-side bias and in the heads-tails bias, we specified $\text{Gamma}(4, 200)$ prior distributions. This choice was motivated by the following considerations: (a) under the $\text{Gamma}(4, 200)$ priors, the expected heterogeneity of the biases (i.e., the expected standard deviation) equals 0.005 when transformed to the probability scale--- half of the hypothesized effect; and (b) under the $\text{Gamma}(4, 200)$ priors, the standard deviation of the heterogeneity equals 0.0025 when transformed to the probability scale--- half of the expected heterogeneity itself. Although these gamma prior distributions were not directly informed by the DHM theory, they are consistent with the estimates from Diaconis et al. \cite{diaconis2007dynamical}; we consider these prior distributions reasonable and relevant for a hypothesis test.

\subsubsection*{Bayesian model-averaging}
The model from Equation~\ref{eq:model} was used to estimate the same-side bias and the heads-tails bias while taking into account the heterogeneity between people and coins. Furthermore, we used the model as a starting point for testing the hypotheses of the same-side bias, $\mathcal{H}_{\beta}$, and heads-tails bias, $\mathcal{H}_{\alpha}$, with the addition of hypotheses about between-people heterogeneity in the same-side bias, $\mathcal{H}_{\sigma_\theta}$, and between-coin heterogeneity in the heads-tails bias, $\mathcal{H}_{\sigma_\gamma}$:
\begin{equation}
    \label{eq:priors-test}
    \begin{split}
        & \mathcal{H}_{\beta_\mu, 1}        : \beta_\mu \sim \text{Beta}(5100, 4900)_{[0.5, 1]} \\
        & \mathcal{H}_{\alpha_\mu, 1}       : \alpha_\mu \sim \text{Beta}(5000, 5000)\\
        & \mathcal{H}_{\sigma_\beta, 1}     : \sigma_\beta \sim \text{Gamma}(4, 200) \\
        & \mathcal{H}_{\sigma_\alpha, 1}    : \sigma_\alpha \sim \text{Gamma}(4, 200) 
    \end{split}
    \begin{split}
        \,\, \text{vs.} \,\, \\
        \,\, \text{vs.} \,\, \\
        \,\, \text{vs.} \,\, \\
        \,\, \text{vs.} \,\,
    \end{split}
    \begin{split}
         & \mathcal{H}_{\beta_\mu, 0}      : \beta_\mu = 0.5 \\
         & \mathcal{H}_{\alpha_\mu, 0}     : \alpha_\mu = 0.5 \\
         & \mathcal{H}_{\sigma_\beta, 0}   : \sigma_\beta = 0 \\
         & \mathcal{H}_{\sigma_\alpha, 0}  : \sigma_\theta = 0. \\ 
    \end{split}
\end{equation}

To test the hypotheses while accounting for uncertainty in the model structure, we used Bayesian model averaging \citep{raftery1995accounting, hoeting1999bayesian,hinne2019conceptual} and specified 16 possible models as an orthogonal combination of the different possible hypotheses. For example, $\mathcal{M}_{1}$ specifies the presence of the same-side bias, the presence of the heads-tails bias, the presence of between-people heterogeneity in same-side bias, and the presence of between-people heterogeneity in the heads/tails bias. $\mathcal{M}_{2}$ specifies the presence of the same-side bias, the presence of the heads-tails bias, the presence of between-people heterogeneity in same-side bias, and the absence of between-people heterogeneity in the heads-tails bias. The last model, $\mathcal{M}_{16}$, then specifies the absence of the same-side bias, the absence of the heads-tails bias, the absence of between-people heterogeneity in same-side bias, and the absence of between-people heterogeneity in the heads/tails bias. The entire model space is listed as follows: 
\begin{align*}
    \mathcal{M}_{1}  &= \mathcal{H}_{\beta_\mu, 1}, \mathcal{H}_{\alpha_\mu, 1},   \mathcal{H}_{\sigma_\beta, 1} \text{, and } \mathcal{H}_{\sigma_\alpha, 1} \\ \nonumber
    \mathcal{M}_{2}  &= \mathcal{H}_{\beta_\mu, 1}, \mathcal{H}_{\alpha_\mu, 1},   \mathcal{H}_{\sigma_\beta, 1} \text{, and } \mathcal{H}_{\sigma_\alpha, 0} \\ \nonumber
    \mathcal{M}_{3}  &= \mathcal{H}_{\beta_\mu, 1}, \mathcal{H}_{\alpha_\mu, 1},   \mathcal{H}_{\sigma_\beta, 0} \text{, and } \mathcal{H}_{\sigma_\alpha, 1} \\ \nonumber
    \mathcal{M}_{4}  &= \mathcal{H}_{\beta_\mu, 1}, \mathcal{H}_{\alpha_\mu, 1},   \mathcal{H}_{\sigma_\beta, 0} \text{, and } \mathcal{H}_{\sigma_\alpha, 0} \\ \nonumber
    \mathcal{M}_{5}  &= \mathcal{H}_{\beta_\mu, 1}, \mathcal{H}_{\alpha_\mu, 0},   \mathcal{H}_{\sigma_\beta, 1} \text{, and } \mathcal{H}_{\sigma_\alpha, 1} \\ \nonumber
    \mathcal{M}_{6}  &= \mathcal{H}_{\beta_\mu, 1}, \mathcal{H}_{\alpha_\mu, 0},   \mathcal{H}_{\sigma_\beta, 1} \text{, and } \mathcal{H}_{\sigma_\alpha, 0} \\ \nonumber
    \mathcal{M}_{7}  &= \mathcal{H}_{\beta_\mu, 1}, \mathcal{H}_{\alpha_\mu, 0},   \mathcal{H}_{\sigma_\beta, 0} \text{, and } \mathcal{H}_{\sigma_\alpha, 1} \\ \nonumber
    \mathcal{M}_{8}  &= \mathcal{H}_{\beta_\mu, 1}, \mathcal{H}_{\alpha_\mu, 0},   \mathcal{H}_{\sigma_\beta, 0} \text{, and } \mathcal{H}_{\sigma_\alpha, 0} \\ \nonumber
    \mathcal{M}_{9}  &= \mathcal{H}_{\beta_\mu, 0}, \mathcal{H}_{\alpha_\mu, 1},   \mathcal{H}_{\sigma_\beta, 1} \text{, and } \mathcal{H}_{\sigma_\alpha, 1} \\ \nonumber
    \mathcal{M}_{10} &= \mathcal{H}_{\beta_\mu, 0}, \mathcal{H}_{\alpha_\mu, 1},   \mathcal{H}_{\sigma_\beta, 1} \text{, and } \mathcal{H}_{\sigma_\alpha, 0} \\ \nonumber
    \mathcal{M}_{11} &= \mathcal{H}_{\beta_\mu, 0}, \mathcal{H}_{\alpha_\mu, 1},   \mathcal{H}_{\sigma_\beta, 0} \text{, and } \mathcal{H}_{\sigma_\alpha, 1} \\ \nonumber
    \mathcal{M}_{12} &= \mathcal{H}_{\beta_\mu, 0}, \mathcal{H}_{\alpha_\mu, 1},   \mathcal{H}_{\sigma_\beta, 0} \text{, and } \mathcal{H}_{\sigma_\alpha, 0} \\ \nonumber
    \mathcal{M}_{13} &= \mathcal{H}_{\beta_\mu, 0}, \mathcal{H}_{\alpha_\mu, 0},   \mathcal{H}_{\sigma_\beta, 1} \text{, and } \mathcal{H}_{\sigma_\alpha, 1} \\ \nonumber
    \mathcal{M}_{14} &= \mathcal{H}_{\beta_\mu, 0}, \mathcal{H}_{\alpha_\mu, 0},   \mathcal{H}_{\sigma_\beta, 1} \text{, and } \mathcal{H}_{\sigma_\alpha, 0} \\ \nonumber
    \mathcal{M}_{15} &= \mathcal{H}_{\beta_\mu, 0}, \mathcal{H}_{\alpha_\mu, 0},   \mathcal{H}_{\sigma_\beta, 0} \text{, and } \mathcal{H}_{\sigma_\alpha, 1} \\ \nonumber
    \mathcal{M}_{16} &= \mathcal{H}_{\beta_\mu, 0}, \mathcal{H}_{\alpha_\mu, 0},   \mathcal{H}_{\sigma_\beta, 0} \text{, and } \mathcal{H}_{\sigma_\alpha, 0}.
\end{align*}

\subsection*{Practice effects}
For both the initial probability of the same-side we side $\beta_{\mu1}$ we specified $\text{Beta}(312, 312)$ prior distribution as in parameter estimation tailored model accounting for heterogeneity. For the between-people heterogeneity in the same-side bias, $\sigma_\beta$, we specified $\text{Normal}_{+}(0, 0.20)$ prior distributions. The $\text{Normal}_{+}(0, 0.20)$ prior distribution results in approximately half-Normal prior distributions with a standard deviation of 0.05 when transformed to probability scale; we chose this values as we wanted to allow a bit more flexibility in the heterogeneity estimate in the practice effect model. For the overall drift, $\theta_\mu$, and the between-people heterogeneity in the overall drift, $\sigma_\theta$, we specified a $\text{Normal}(0, 0.10)$ and  $\text{Normal}_{+}(0, 0.10)$ prior distributions which allowed for enough flexibility when estimating the per 100 flips aggregated drift. For the overall state change at flip i, $\delta_{\mu i }$, and the between-people heterogeneity in the state change at flip i, $\sigma_\delta$, we specified a $\text{Normal}(0, 0.20)$ and  $\text{Normal}_{+}(0, 0.20)$ prior distributions which again allowed for enough flexibility when estimating the per 100 flips state change. Later visualizations confirmed that the specified prior distributions were not overly restrictive.

\subsection*{Computational details}
All mixed-effect models were estimated in \texttt{Stan} \citep{Stan} via the \texttt{rstan} \texttt{R} package \citep[version 2.26.1][]{Rstan}. The marginal likelihood was determined using bridge sampling \citep{meng1996simulating, gronau2017tutorial} via the \texttt{bridgesampling} \texttt{R} package \citep[version 1.1-2]{bridgesampling}. The mixed-effect models were run using $10$ chains with $15{,}000$ warm-up iterations and $10{,}000$ sampling iterations each. All models converged with $\Hat{R} < 1.01$. Only the state-space model showed a couple of divergent transitions.

\section*{Appendix B: Prior sensitivity analysis}
\label{app:prior_sensitivity}

We assessed the sensitivity of the performed hypothesis tests to the prior specification. First, we examined prior sensitivity of the Bayesian (nonhierarchical) informed binomial hypothesis tests; second, we examined prior sensitivity of the more complex hierarchical Bayesian logistic regression utilizing Bayesian model averaging. We specified a range of prior distributions for each alternative hypothesis by employing the normal-moment prior distribution \citep{johnson2010use}, Normal-Moment($\phi$),
\begin{equation*}
    p(x \mid \phi) = \frac{2x^2}{ \sqrt{\pi}  |\phi^3|} \text{exp}\left(\frac{-x^2}{\phi^2} \right),
\end{equation*}
which allows us to specify non-local prior distributions using a single parameter---the mode of the prior distribution $\phi$. When the support of the Normal-Moment distribution is unrestricted, it features two modes at $\pm \phi$ with the density decreasing towards zero at 0 and towards zero at $\pm \infty$. Varying the $\phi$ parameter across a range of plausible values allows us to assess the evidence in favor or against the hypothesis as a function of the prior distribution's mode, sometimes described as the Bayes factor function \citep[BFF,][]{johnson2023bayes}. We examined the prior sensitivity with both the full data set and on data set after excluding potential outliers (i.e., people with the observed probability of the same side larger than 0.53). 

All prior sensitivity analyses were conducted using \texttt{Stan} \citep{Stan} via the \texttt{rstan} \texttt{R} package \citep[version 2.26.1][]{Rstan} and estimated the marginal likelihood using bridge sampling \citep{meng1996simulating, gronau2017tutorial} using the \texttt{bridgesampling} \texttt{R} package \citep[version 1.1-2]{bridgesampling}. The models were run using $4$ chains with $15{,}000$ warm-up and $10{,}000$ sampling iterations.

For an alternative approach (an inversion of Bayes factor function) see Section~4.1 and 4.2 in Pawel \cite{pawel2024bayes}, which leads to similar conclusions.

\subsection*{Nonhierarchical analyses}

Using the Normal-Moment prior distributions results in the modified definition of our hypotheses about the same-side bias (Equation~\ref{eq:overal_same_side}):
\begin{align*}
    \text{No same-side bias, } \mathcal{H}_{0}    &: \text{logit}(\beta) = 0.5 \\ \nonumber
    \text{DHM same-side bias, } \mathcal{H}_{1} &: \text{logit}(\beta) \sim\text{Normal-Moment}_+(\phi_\beta), \nonumber
\end{align*}
and the heads-tails bias:
\begin{align*}
    \text{No heads-tails bias, } \mathcal{H}_{0a}       &: \text{logit}(\alpha) = 0.5 \\ \nonumber
    \text{Small heads-tails bias, } \mathcal{H}_{1a}    &: \text{logit}(\alpha) \sim\text{Normal-Moment}(\phi_\alpha). \nonumber
\end{align*}
The main difference lies in defining prior distributions for the same-side bias $\text{logit}(\beta)$ and heads-tails bias $\text{logit}(\alpha)$ instead of the probability of same side $\beta$ and probability of heads $\alpha$.

We examined a range of modes $\phi_\beta = [0, 0.08]$ and $\phi_\alpha = [0, 0.08]$ of the Normal-Moment distributions. The examined modes translate to the mode probability of the same side ranging from $0.50$ to $0.52$ and the probability of heads ranging from $0.48$ to $0.52$. The upper ranges of the prior modes were chosen a posteriori based on the observed ranges of estimates 

\subsubsection*{Results}

\begin{figure}[h]
    \centering
    \includegraphics[width=1\textwidth]{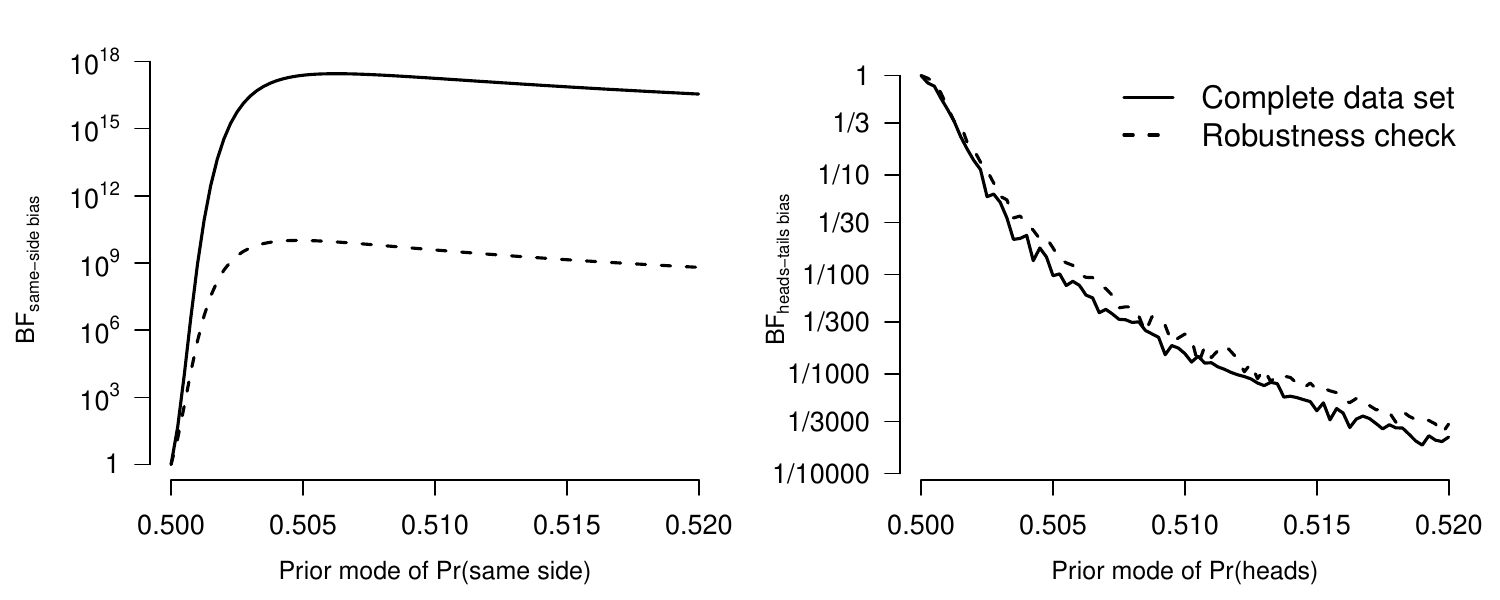}
    \caption{A prior sensitivity analysis demonstrates robustness of the same-side bias (left) and the absence of the heads-tails bias (right) across different prior specifications ($x$-axis) and exclusion of potential outliers (full line vs. dashed line) from nonhierarchical models.}
    \label{fig:robustness1}
\end{figure}

The left panel of Figure~\ref{fig:robustness1} visualizes the Bayes factor functions for the same-side bias using the complete data set (full line) and after excluding the potential outliers (dashed line). The evidence in favor of the same-side bias is rapidly increasing even for small deviations from the null hypothesis; we find strong evidence ($\text{BF}_{\text{same-side bias}} > 10$) for all examined modes of alternative hypothesis of the presence of the same-side bias (i.e., $\text{Pr}(\text{same side} > 0.5003)$) regardless of potential outlier exclusion. The alternative hypothesis of the same-side bias receives the highest support at the mode $\text{Pr}(\text{same side} = 0.5062)$ ($\text{BF}_{\text{same-side bias}} = 2.85 \times 10^{17}$) using the full data set and at the mode $\text{Pr}(\text{same side} = 0.5047)$ ($\text{BF}_{\text{same-side bias}} = 1.01 \times 10^{10}$) after excluding potential outliers.

The right panel of Figure~\ref{fig:robustness1} visualizes the Bayes factor functions for the heads-tails bias using the complete data set (full line) and after excluding the potential outliers (dashed line). The evidence against the heads-tails bias is increasing for any deviations from the null hypothesis; we find strong evidence against the heads-tails bias ($\text{BF}_{\text{heads-tails bias}} < 1/10$) for all examined modes of alternative hypothesis of the heads-tails bias starting at $\text{Pr}(\text{heads} > 0.5025)$ using the full data set and at $\text{Pr}(\text{same side} > 0.5028)$ after excluding potential outliers.

\subsection*{Analyses accounting for dependencies due to people and coins}

Using the Normal-Moment prior distributions results in the modified definition of our hypotheses (Equation~\ref{eq:priors-test}):
\begin{equation}
    \label{eq:priors-test-robustness}
    \begin{split}
        & \mathcal{H}_{\beta_\mu, 1}        : \text{logit}(\beta_\mu)  \sim \text{Normal-Moment}_+(\phi_{\beta_\mu})           \\
        & \mathcal{H}_{\alpha_\mu, 1}       : \text{logit}(\alpha_\mu) \sim \text{Normal-Moment}(\phi_{\alpha_\mu})            \\
        & \mathcal{H}_{\sigma_\theta, 1}: \sigma_\theta        \sim \text{Normal-Moment}_+(\phi_{\sigma_\theta}) \\
        & \mathcal{H}_{\sigma_\gamma, 1}: \sigma_\gamma        \sim \text{Normal-Moment}_+(\phi_{\sigma_\gamma}) 
    \end{split}
    \begin{split}
        \,\, \text{vs.} \,\, \\
        \,\, \text{vs.} \,\, \\
        \,\, \text{vs.} \,\, \\
        \,\, \text{vs.} \,\,
    \end{split}
    \begin{split}
         & \mathcal{H}_{\beta_\mu, 0}        : \beta_\mu = 0.5 \\
         & \mathcal{H}_{\alpha_\mu, 0}       : \alpha_\mu = 0.5 \\
         & \mathcal{H}_{\sigma_\theta, 0}: \sigma_\theta = 0 \\
         & \mathcal{H}_{\sigma_\gamma, 0}: \sigma_\gamma = 0. \\ 
    \end{split}
\end{equation}
We examined the sensitivity of each hypothesis test marginally to the non-assessed hypotheses (i.e., we fixed the prior mode of all but the tested parameters) to keep the computational load feasible. This corresponds to computing four BFFs under the corresponding settings for the alternative hypotheses (the null hypotheses remain unchanged): 
\begin{align*}
    \text{BFF}_{\beta_\mu}  &: \phi_{\beta_\mu} = [0, 0.08], \phi_{\alpha_\mu} = 0.04, \phi_{\sigma_\theta} = 0.02, \text{and } \phi_{\sigma_\gamma} = 0.02     \\ \nonumber
    \text{BFF}_{\alpha_\mu} &: \phi_{\beta_\mu} = 0.04, \phi_{\alpha_\mu} = [0, 0.08], \phi_{\sigma_\theta} = 0.02, \text{and } \phi_{\sigma_\gamma} = 0.02     \\ \nonumber
    \text{BFF}_{\sigma_\beta}  &: \phi_{\beta_\mu} = 0.04, \phi_{\alpha_\mu} = 0.04, \phi_{\sigma_\theta} = [0, 0.08], \text{and } \phi_{\sigma_\gamma} = 0.02     \\ \nonumber
    \text{BFF}_{\sigma_\beta}  &: \phi_{\beta_\mu} = 0.04, \phi_{\alpha_\mu} = 0.04, \phi_{\sigma_\theta} = 0.02, \text{and } \phi_{\sigma_\gamma} = [0, 0.08].     
\end{align*}
The fixed prior modes were chosen to correspond approximately to the mean of prior distributions specified in Equation~\ref{eq:priors-test}, that is, a same-side bias $\beta \approx 0.51$, heads-tails bias $\alpha \approx 0.49 \text{ and } 0.51$, between-people heterogeneity in the same-side bias $\sigma_\theta \approx 0.01$, and between-coin heterogeneity in the heads-tails bias $\sigma_\theta \approx 0.01$. The upper ranges of the prior modes were chosen a posteriori, based on the observed ranges of estimates, i.e., up to a bias of $\approx 0.52$ or between-people / coin heterogeneity in the biases of $\approx 0.02$.

\subsubsection*{Results}

\begin{figure}[h!]
    \centering
    \includegraphics[width=1\textwidth]{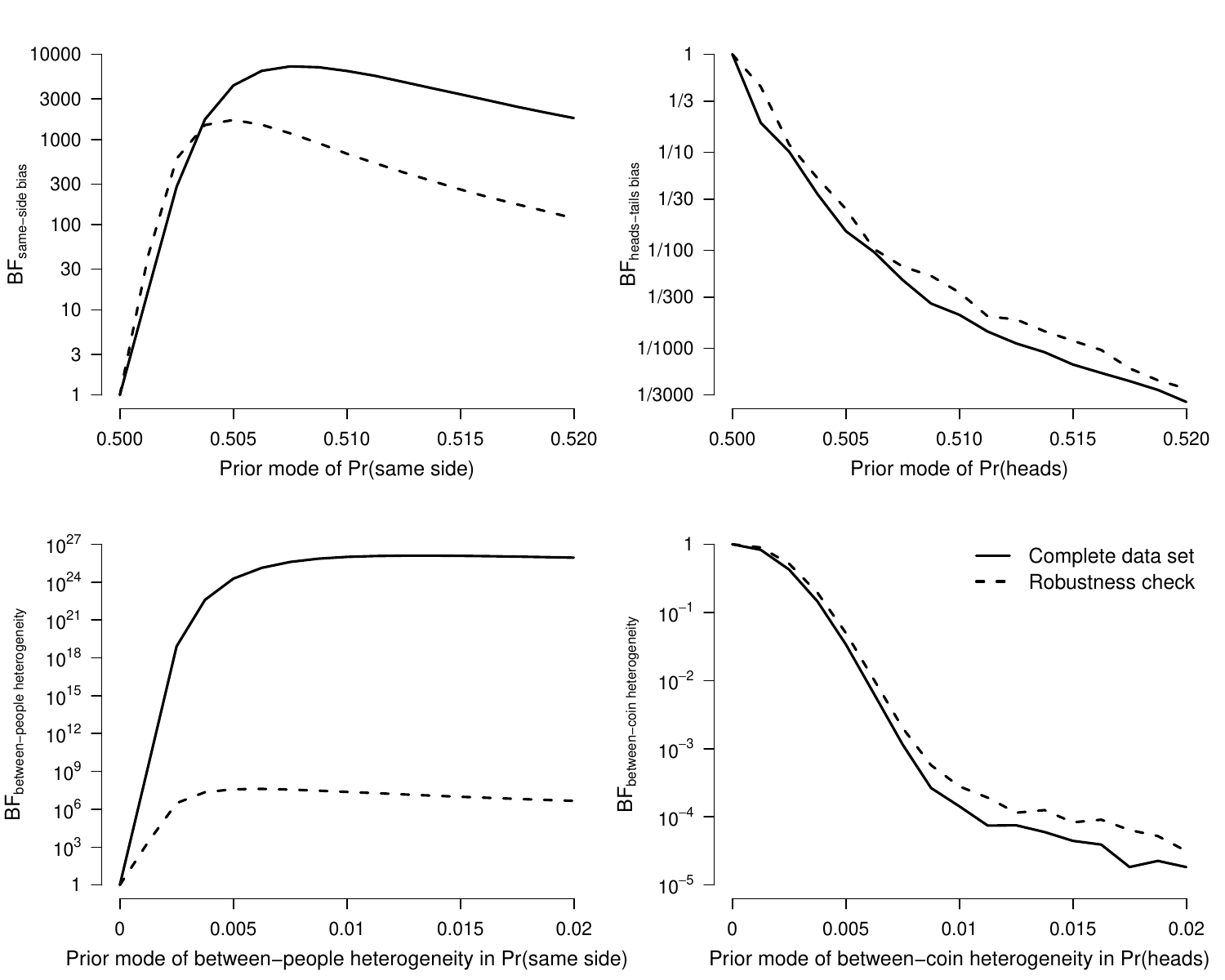}
    \caption{Prior sensitivity analysis demonstrates robustness of the same-side bias (top left), the absence of the heads-tails bias (top right), the presence of between-people heterogeneity in the same-side bias (bottom left), and the absence of between-coin heterogeneity in the heads-tails bias (bottom left) across different prior specifications ($x$-axis) and exclusion of the potential outliers (full vs. dashed line) from hierarchical models accounting for heterogeneity between people and coins.}
    \label{fig:robustness2}
\end{figure}

The first row of Figure~\ref{fig:robustness2} visualizes the Bayes factor functions for the same-side bias (left) and heads-tails bias (right) when accounting for between-people and between-coin heterogeneity using the complete data set (full line) and after excluding the potential outliers (dashed line). In contrast to the heterogeneity unadjusted analysis (cf., Figure~\ref{fig:robustness1}), the degree of evidence in favor of the same-side bias and against heads-tails bias is notably smaller. Nevertheless, there is strong evidence ($\text{BF}_{\text{same-side bias}} > 10$) for all examined modes of alternative hypothesis of the presence of the same-side bias (i.e., $\text{Pr}(\text{same side} > 0.5013)$) regardless of potential outlier exclusion. The alternative hypothesis of the same-side bias receives the highest support at the mode $\text{Pr}(\text{same side} = 0.5075)$ ($\text{BF}_{\text{same-side bias}} = 7208$) using the full data set and at the mode $\text{Pr}(\text{same side} = 0.5050)$ ($\text{BF}_{\text{same-side bias}} = 1697$) after excluding potential outliers. The heads-tails bias is decreasing slower than in the heterogeneity unadjusted analysis and we find strong evidence for the absence of heads-tails bias ($\text{BF}_{\text{heads-tails bias}} < 1/10$) only against alternative hypotheses with mode $\text{Pr}(\text{heads} = 0.5037)$ or larger regardless of potential outlier exclusion.

The second row of Figure~\ref{fig:robustness2} visualizes the Bayes factor functions for the between-people heterogeneity in the same-side bias (left) and the between-coin heterogeneity in heads-tails bias (right) using the complete data set (full line) and after excluding the four potential outliers (dashed line). The degree of evidence in favor of the between-people heterogeneity in the same-side bias is extreme: we find strong evidence for the between-people heterogeneity in the same-side bias in all examined modes of alternative hypothesis, with the strongest evidence for the alternative hypothesis of the between-people heterogeneity at the mode $\text{sd}_\text{people}(\text{Pr}(\text{same side})) = 0.0137$ ($\text{BF}_{\text{between-people heterogeneity}} = 1.27 \times 10^{26}$) using the full data set and at the mode $\text{sd}_\text{people}(\text{Pr}(\text{same side})) = 0.0062$ ($\text{BF}_{\text{between-people heterogeneity}} = 4.15 \times 10^{7}$) after excluding potential outliers. The evidence for the absence of between-coin heterogeneity in heads-tails bias accumulates more slowly: we find strong evidence for the absence of between-coin heterogeneity in heads-tails bias ($\text{BF}_{\text{between-coin heterogeneity}} < 1/10$) only against alternative hypotheses with mode $\text{sd}_\text{coins}(p_{\text{heads}}) = 0.0050$ or larger regardless of potential outlier exclusion.

\subsection*{Summary}

The prior sensitivity analysis verified that our findings are robust to the specific choices of our prior distributions. In fact, the prior sensitivity analysis suggests that the data support a substantive range of parameterization of the same-side bias and between-people heterogeneity in the same-side bias.

\section*{Appendix C: Same-side bias variability by recruitment site}

We assessed whether the degree of the same-side bias systematically varies according to the source of participant's recruitment. We extended the same-side bias part of the Equation~\ref{eq:model}, $\text{logit}(\beta)_k$, with a by-recruitment site difference from the same-side bias $\delta_{k}$,
\begin{equation*}
    \label{eq:model-reg}
    \text{logit}(\beta_k)   =    \text{logit}(\beta_\mu) + \gamma_{\beta_k} + \delta_{k},    \\  \nonumber
\end{equation*}
such that the $\delta_{k}$ implement a sum to zero constrain via standardized orthonormal contrast.

We only attempted to estimate the parameters and specified slightly informed prior distributions for $\alpha_\mu$, $\beta_\mu$, $\sigma_\alpha$, and $\sigma_\beta$ and a slightly informed prior distribution for the by-recruitment site difference from the same-side bias  
\begin{equation*}
    \delta_{k}      \sim  \text{Normal}(0, 0.20),
\end{equation*}
which results in standard deviation of $0.05$ on probability scale. We chose this prior distribution over our originally intended prior distribution $\text{Normal}(0, 0.04)$ (i.e., leading to a $0.01$ standard deviation on the probability scale), because the originally intended distribution turned out to restrict the posterior distribution too severely (the data contained very little information about the by-recruitment site differences).

\subsection*{Results}

\begin{figure}[h]
    \centering
    \includegraphics[width=0.75\textwidth]{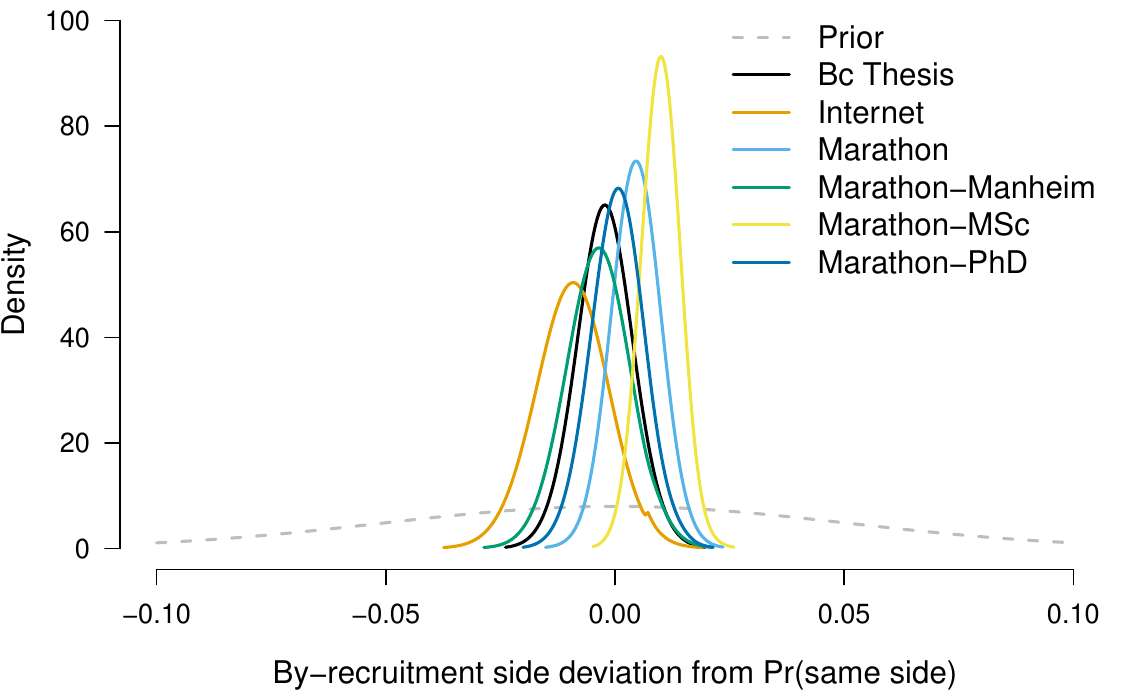}
    \caption{Posterior distribution of by-recruitment site differences from the probability of the same side show a considerable amount of uncertainty but no systematic differences. Prior distribution as grey dashed line, by-recruitment side posterior distributions in colors. Only the 95\% CI for `Marathon-MSc' does not contain zero; however, the difference disappears after excluding outliers, see Online Supplement Figure~3.}
    \label{fig:regression}
\end{figure}

Figure~\ref{fig:regression} visualizes the prior and posterior distributions of by-recruitment site differences from the probability of the same side. Note that all posterior distributions remain relatively wide, highlighting a lack of information about potential person-level moderators. Out of all six sites, only the 95\% CI for the `Marathon-Msc' does not contain zero, $\delta_\text{Marathon-Msc} = 0.0010$, 95\% CI [$0.0014$, $0.0184$], however, the difference disappears after excluding outliers. The remaining by-recruitment site differences: $\delta_\text{Bc Thesis} = -0.0022$, 95\% CI [$-0.0146$, $0.0103$], $\delta_\text{Internet} = -0.0093$, 95\% CI [$-0.0253$, $0.0066$], $\delta_\text{Marathon} = 0.0046$, 95\% CI [$-0.0063$, $0.0153$], $\delta_\text{Marathon-Manheim} = -0.0036$, 95\% CI [$-0.0177$, $0.0105$], and $\delta_\text{Marathon-PhD} = 0.0006$, 95\% CI [$-00112.$, $0.0124$].

Removal of the four potential outliers (i.e., people with an observed probability of the same side larger than 0.53) resulted in a considerable decrease in the by-recruitment site differences from the probability of the same side. The `Marathon-Msc`'s 95\% CI of by-recruitment site differences from the probability of the same side now contains zero, $\delta_\text{Marathon-Msc} = 0.0016$, 95\% CI [$-0.0037$, $0.0069$], however, the by recruitment site difference for `Marathon` no longer includes zero $\delta_\text{Marathon} = 0.0062$, 95\% CI [$0.0004$, $0.0120$]. See Figure~4 in Online Supplements for the corresponding visualization.

\subsection*{Summary}

The analysis did not reveal reliable and systematic by-recruitment site differences from the same-side bias; `Marathon-Msc' site deviation from the same-side bias's 95\% CI did not contain zero using the complete dataset, and `Marathon' site deviation from the same-side bias's 95\% CI did not contain zero using the potential outliers removed data set. The posterior estimates were however considerably variable and a much larger data set would be needed to draw strong conclusions.

\section*{Appendix D: Frequentist results}
\label{app:frequentist}

At the request of editors and reviewers, we reproduce the main findings from our manuscript using a frequentist methodology. All CIs reported in this section correspond to 95\% confidence intervals. In each subsection, we report results from an analysis performed on the full data set and follow-up results performed on the data set after excluding potential outliers (i.e., the four people with an observed probability of the same side larger than 0.53). All mixed-effects models were estimated using the \texttt{lme4} \texttt{R} package \citep[version 1.1.35.3,][]{lme4}.

\subsection*{Analyses assuming independent observations}

Assuming independence of the individual coin tosses, a simple (nonhierarchical) binomial test finds a statistically significant effect of the same-side bias ($178{,}079$ same-side landings from $350{,}757$ tosses), $p < 0.001$, $\text{Pr}(\text{same side}) = 0.5077$, 95\% CI [$0.5060$, $0.5094$]. Furthermore, a simple binomial test does not find a statistically significant effect of a heads-tails bias ($175{,}421$ heads out of $350{,}757$ tosses), $p = 0.887$, $\text{Pr}(\text{heads}) = 0.5001$, 95\% CI [$0.4985$, $0.5018$].

Removing potential outliers results in a same-side bias that is smaller, $\text{Pr}(\text{same side}) = 0.5060$, 95\% CI [$0.5043$, $0.5077$] ($171{,}517$ same-side landings from $338{,}985$ tosses), but remains statistically significant, $p < 0.001$. The probability of heads remains practically identical, $\text{Pr}(\text{heads}) = 0.5004$, 95\% CI [$0.4987$, $0.5021$] ($169{,}635$ heads out of $338{,}985$ tosses), and statistically non-significant, $p = 0.626$.

\subsection*{Analyses accounting for dependencies due to people and coins}

To account for the possible heterogeneity between coins and people, we specified a generalized linear mixed-effect model with a binomial likelihood and a logit link for the probability of a coin landing on the same side it started with: the intercept ($b_\mu$, corresponding to the same-side bias), by-participant random intercepts ($\tau_{b_\mu}$, corresponding to the by-participant variability in the same-side bias), and the fixed-effect of starting side ($b_1$, a factor with the following levels: `heads', `tails'). The model reveals a statistically significant same-side bias $b_\mu = 0.0399$, $\text{se} = 0.0099$, $z = 4.03$, $p < 0.001$ (corresponding to the probability of the same side $\text{Pr}(\text{same side}) = 0.5100$, 95\% CI [$0.5051$, $0.5148$]) and a statistically non-significant effect of the starting side, $b_1 = 0.0007$, $\text{se} = 0.0034$, $z = 0.22$, $p = 0.829$. The model also reveals a statistically significant between-people variability in the same-side bias, $\tau_{b_\mu} = 0.0626$, $\chi^2(1) = 120.43$, $p < 0.001$ (by a comparison to a model without the random intercept). We could not test for a between-coin variance in the heads-tails bias due to the limited flexibility of \text{lme4}. However, we tested for a potential between-coin variance in the same-side bias---the model extension with the by-coin intercept in the same-side bias turns out statistically non-significant, $\chi^2(1) = 1.19$, $p = 0.276$.

Removing potential outliers does not qualitatively affect the results. The degree of the same-side bias decreases but remains statistically significant $b_\mu = 0.0241$, $\text{se} = 0.0057$, $z = 4.24$, $p < 0.001$  (corresponding to the probability of the same side $\text{Pr}(\text{same side}) = 0.5060$, 95\% CI [$0.5032$, $0.5088$]) and the effect of the starting side remains statistically non-significant, $b_1 = 0.0018$, $\text{se} = 0.0034$, $z = 0.52$, $p = 0.604$. The model again reveals a smaller but statistically significant between-people variability in the same-side bias, $\tau_{b_\mu} = 0.0278$, $\chi^2(1) = 36.05$, $p < 0.001$ (by a comparison to a model without the random intercept), and statistically non-significant extension with a between-coin variance in the same-side bias, $\chi^2(1) = 2.77$, $p = 0.096$.

We further refer the reader to the \texttt{dat.bartos2023} entry in the \texttt{metadat} \citep[version 1.3-0,][]{metadat} \texttt{R} package by Viechtbauer \cite{viechtbauerVignette}, who re-analyzed the data set with a frequentist meta-analyses using the \texttt{metafor} \citep{metafor} \texttt{R} package and arrived at qualitatively identical conclusions.  

\subsection*{Practice effects}

\begin{figure}[h]
    \centering
    \includegraphics[width=1\textwidth]{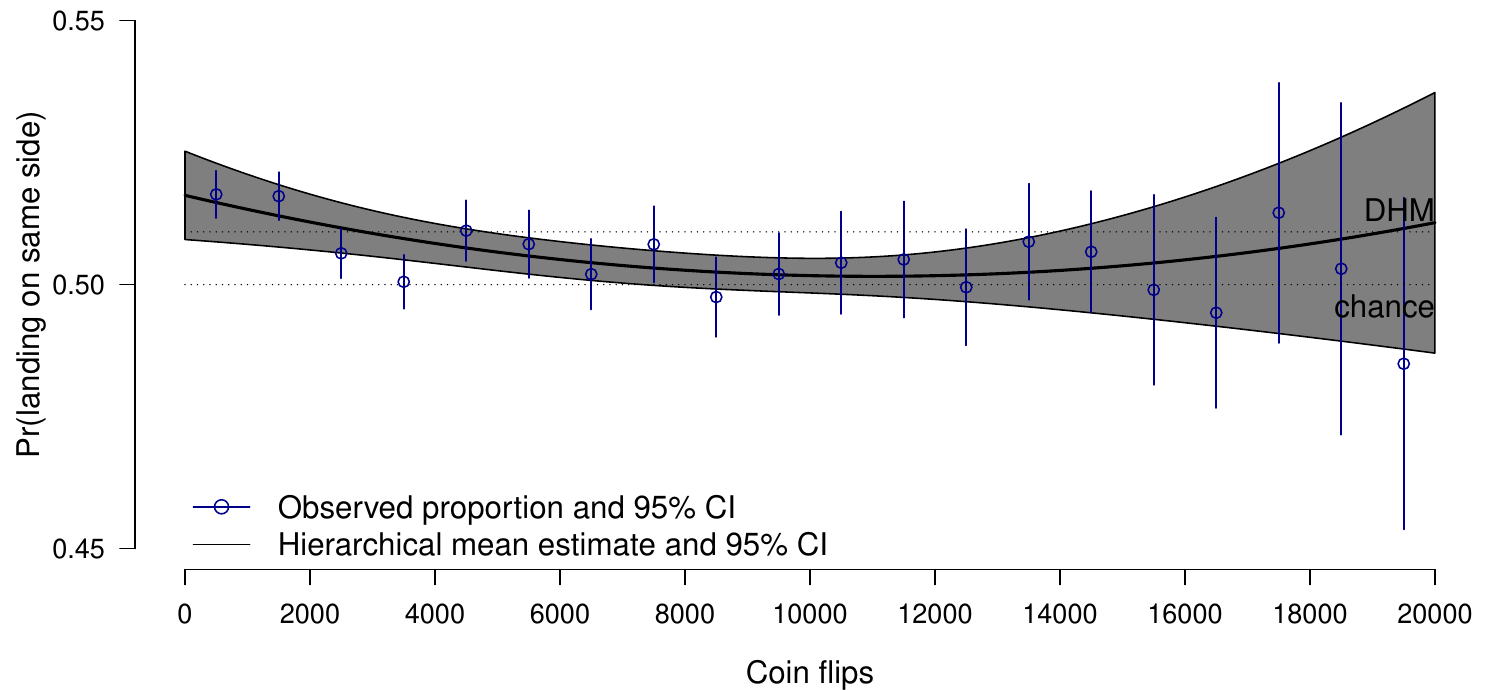}
    \caption{Frequentist analysis verifying the degree of the individual same-side bias decreases as people flip more coins. Mean estimate and 95\% CI of the probability of the same side with the increasing number of coin flips from a generalized linear mixed effect model with binomial likelihood, logit link function, and linear and quadratic effect of time on the probability of the same side, vs. the observed proportion and 95\% CI of the same side probability aggregated per sequences of $1{,}000$ coin flips.}
    \label{fig:learning-f}
\end{figure}

To estimate the learning effects, we extended the previously specified generalized linear mixed-effect model by including the linear ($\beta_t$) and quadratic ($\beta_{t^2}$ fixed-effect of the number of coin flips (scaled by $10{,}000$ to improve convergence), and the corresponding by-participant random slopes ($\tau_{b_t}$ and $\tau_{b_{t^2}}$; and the corresponding correlation between random-effects). This model allows us to estimate and test the initial same-side bias (the intercept $b_\mu$) but the polynomial form does not provide an estimate of the asymptotic same-side bias. We first tested linear fixed- and random-effects of time and found that the model fit significantly improves upon the previously specified linear mixed-effect model, $\chi^2(3) = 45.42$, $p < 0.001$. Then we tested additional quadratic fixed- and random-effects of time and found the model fit significantly improves upon the linear fixed- and random-effects only model, $\chi^2(4) = 23.47$, $p = 0.001$. We did not pursue more complex models and report the linear and quadratic fixed and random-effects model. The model reveals a statistically significant initial same-side bias much larger than the previous models, $b_\mu = 0.0677$, $\text{se} = 0.0171$, $z = 3.96$, $p < 0.001$ and substantial between-people variability in the initial same-side bias $\tau_{b_\mu} = 0.1048$. The statistically significant linear fixed-effect of time, $b_t = -0.1114$, $\text{se} = 0.0454$, $z = -2.45$, $p = 0.014$, and statistically non-significant quadratic fixed-effect of time, $b_{t^2} = 0.0505$, $\text{se} = 0.0305$, $z = 1.66$, $p = 0.098$, are both accompanied by notably between-people variability in the corresponding effects, $\tau_{b_t} = 0.2196$ and $\tau_{b_{t^2}} = 0.1241$. The time effects result in a high initially probability of the same side, $\text{Pr}(\text{initial same side}) = 0.5169$, 95\% CI [$0.5085$, $0.5253$], which decreases quickly in the initial couple thousands flips and  eventually starts increasing (with an extreme uncertainty) due to the quadratic form of the time effect (Figure~\ref{fig:learning-f}). The results again confirm the statistically non-significant effect of the starting side, $b_1 = 0.0008$, $\text{se} = 0.0034$, $z = 0.23$, $p = 0.816$.

Removing the four potential outliers does not qualitatively affect the results. The linear fixed- and random-effects of time model significantly improves upon the simple linear mixed-effect model, $\chi^2(3) = 23.39$, $p < 0.001$; the quadratic fixed- and the random-effects of time model significantly improves upon the linear fixed- and random-effects of time model, $\chi^2(4) = 12.54$, $p = 0.014$. The statistically significant initial same-side bias decreases after excluding the potential outliers, $b_\mu = 0.0443$, $\text{se} = 0.0131$, $z = 3.37$, $p < 0.001$ as well as the between-people variability in the initial same-side bias $\tau_{b_\mu} = 0.0683$. The smaller statistically non-significant linear fixed-effect of time, $b_t = -0.0613$, $\text{se} = 0.0387$, $z = -1.59$, $p = 0.113$, and the smaller statistically non-significant quadratic fixed-effect of time, $b_{t^2} = 0.0233$, $\text{se} = 0.0271$, $z = 0.86$, $p = 0.390$, are still accompanied by smaller but notable between-people variability in the corresponding effects, $\tau_{b_t} = 0.1470$ and $\tau_{b_{t^2}} = 0.0877$. However, the time effects still result in smaller yet high initially probability of the same side, $\text{Pr}(\text{initial same side}) = 0.5111$, 95\% CI [$0.5046$, $0.5175$], which still decreases quickly in the initial couple thousands flips and eventually starts increasing (with an extreme uncertainty) due to the quadratic form of the time effect (Figure~5 in Online Supplements). The results again confirm the statistically non-significant effect of the starting side, $b_1 = 0.0019$, $\text{se} = 0.0034$, $z = 0.54$, $p = 0.588$.

\subsection*{Same-side bias variability by recruitment site}

Finally, we tested for the potential variability in the same-side bias according to the source of participant's recruitment. We specified a generalized linear mixed-effect model with a binomial likelihood and a logit link for the probability of a coin landing on the same side it started with: the intercept ($b_\mu$, corresponding to the same-side bias), by-participant random intercepts ($\tau_{b_\mu}$, corresponding to the by-participant variability in the same-side bias), fixed-effects for the source of participants' recruitment ($b_{r1}, b_{r2}, \dots, b_{r5}$; a factor with the following levels: `Bc Thesis', `Marathon', `Internet', `Marathon-Manheim', `Marathon-MSc', `Marathon-PhD'), and the fixed-effect of starting side ($b_1$, a factor with the following levels: `heads', `tails'). A comparison of the full model and a reduced model (omitting the the source of participants' recruitment) did not reveal a significant effect of the source of participants' recruitment, $\chi^2(5) = 6.65$, $p = 0.248$.

\begin{figure}[h]
    \centering
    \includegraphics[width=0.75\textwidth]{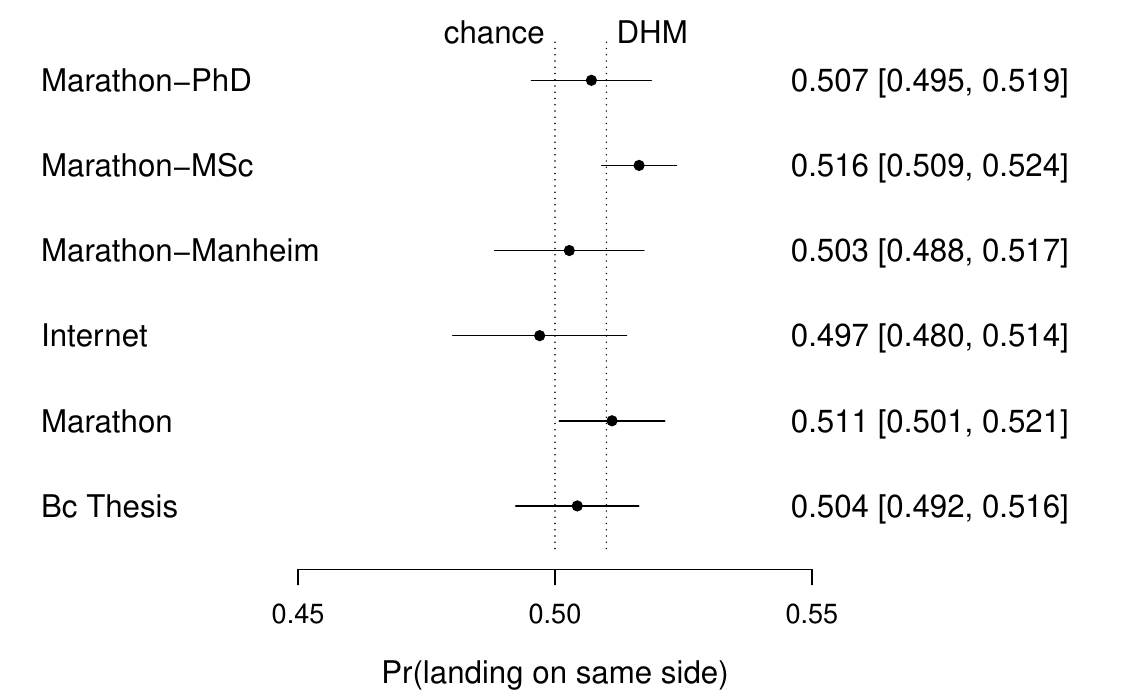}
    \caption{Frequentist marginal mean estimates and 95\% CI of the probability of the same side do not differ by the source of participants' recruitment.}
    \label{fig:recruitment-f}
\end{figure}

Removing potential outliers does not affect the results as the source of recruitment fixed effects remain statistically non-significant, $\chi^2(5) = 7.71$, $p = 0.173$. Figure~6 in Online Supplements visualizes the estimates.

\subsection*{Summary}

All frequentist re-analyses are aligned with the Bayesian results presented in the manuscript; a statistically significant same-side bias, statistically non-significant heads-tails bias (with our without adjustment for between-people heterogeneity), statistically significant between-people heterogeneity in the same-side bias, and statistically significant learning effects that decrease the same-side bias over time. No results change after exclusion of potential outliers.

\section*{Appendix E: Audit}

We randomly sampled and audited 90 sequences of 100 coin flips each. We verified the existence of the video recordings (with occasionally missing video recordings due to file corruption or recording equipment malfunction) and attempted to re-code the outcome of individual coin tosses from the video recordings. We encountered video recordings of varying quality and detail which made one-to-one matching of the original coded sequences and the re-coded audited sequences highly challenging. However, assessing the degree of same-side bias on the original vs. the audited sequences revealed that the original sequences contained a highly similar degree of same-side bias. As such, it seems implausible that the original sequences were affected by coding bias in favor of the same-side hypothesis. 

\end{document}